\documentclass[aos,MSNbibl,seceqn,citesort,dvips]{arximspdf}
\usepackage{dcolumn}
\usepackage{graphicx}

\doi{10.1214/11-AOS936}
\volume{40}
\issue{1}
\pubyear{2012}
\firstpage{1}
\lastpage{29}

\makeatletter

\newcolumntype{d}[1]{D{.}{.}{#1}}

\newcommand{\bs}{\bolds}
\newcommand{\mf}{\mathbf}
\newcommand{\E}{\mathrm{E}}
\newcommand{\Cov}{\operatorname{Cov}}

\newcommand{\R}{\mathbb{R}}

\newcommand{\sumii}{\sum_i}
\newcommand{\sumj}{\sum_{j=2}^{n_i}}
\newcommand{\hutij}{\hat{\mu}(t_{ij})}
\newcommand{\hphik}{\hat{\phi}_k}
\newcommand{\dij}{\tau_{ij}}
\newcommand{\mani}{\mathcal{M}}
\newcommand{\sumk}{\sum_{k\leq K}}
\newcommand{\humani}{\hat{\bolds\mu}}
\newcommand{\umani}{\bolds{\mu}}
\newcommand{\hejmani}{\hat{\mf{e}}_j}
\newcommand{\ejmani}{\mf{e}_j}
\newcommand{\tXit}{\tilde{X}_i(t)}
\newcommand{\hXiK}{\hat{X}_i^K}
\newcommand{\hXmu}{\hat{\mu}^{\mani}}
\newcommand{\Xmu}{\mu^{\mani}}
\newcommand{\hXj}{\hat{X}_{j,\alpha}^{\mani}}
\newcommand{\Xj}{X_{j,\alpha}^{\mani}}
\newcommand{\bu}{\mf{u}}
\newcommand{\slj}{(\lambda_j^{\mani})^{1/2}}

\newcommand{\argmin}[1]{\mathop{\arg\min}_{#1}}

\newtheorem{Thm}{Theorem}

\newtheorem{Cor}{Corollary}
\newtheorem{Prop}{Proposition}

\makeatother

\begin{document}
\begin{frontmatter}

\title{Nonlinear manifold representations for functional~data}
\runtitle{Functional manifold representations}

\begin{aug}
\author[A]{\fnms{Dong} \snm{Chen}\corref{}\ead[label=e1]{dchen@wald.ucdavis.edu}}
\and
\author[A]{\fnms{Hans-Georg} \snm{M\"uller}\thanksref{t1}\ead[label=e2]{mueller@wald.ucdavis.edu}}
\runauthor{D. Chen and H.-G. M\"uller}
\affiliation{University of California, Davis}
\address[A]{Department of Statistics\\
University of California\\
Davis, California 95616\\
USA\\
\printead{e1}\\
\hphantom{E-mail: }\printead*{e2}} %adresu isvedimo komanda gale!
\end{aug}

\thankstext{t1}{Supported in part by NSF Grants DMS-08-06199 and DMS-11-04426.}

% HISTORY:
\received{\smonth{6} \syear{2011}}
\revised{\smonth{9} \syear{2011}}

% ABSTRACT
%
\begin{abstract}
For functional data lying on an unknown nonlinear low-dimen\-sional
space, we study manifold learning and introduce the notions of manifold
mean, manifold modes of functional variation and of functional manifold
components. These constitute nonlinear representations of functional
data that complement classical linear representations such as
eigenfunctions and functional principal components. Our manifold
learning procedures borrow ideas from existing nonlinear dimension
reduction methods, which we modify to address functional data settings.
In simulations and applications, we study examples of functional data
which lie on a manifold and validate the superior behavior of manifold
mean and functional manifold components over traditional
cross-sectional mean and functional principal components. We also
include consistency proofs for our estimators under certain
assumptions.
\end{abstract}

% KEYWORDS
%

\end{frontmatter}

%s1 #&#
\section{Introduction}\label{sec1}

Nonlinear dimension reduction methods, such as locally linear embedding
\cite{row00}, isometric mapping \cite{ten00} and Laplacian eigenmaps
\cite{bel03}, have been successfully applied to image data in recent
years. A commonly used example is the analysis of photos of a sculpture
face taken under different angles and lighting conditions. The number
of pixels of these images is huge, but their structure only depends on
a few variables related to angle and lighting conditions. It is then
advantageous to treat the observed image data as a manifold that is
approximately \textit{isomorphic} to a low-dimensional Euclidean space.

Unlike shape analysis \cite{ken99} and the recent diffusion tensor
imaging \cite{huc11}, where it is assumed that the form of the manifold
is known a priori, nonlinear dimension reduction methods usually are
manifold-learning procedures, where the manifold\vadjust{\goodbreak} is not known but it is
assumed that it possesses certain features which are preserved in the
observed data. For instance, locally linear embedding preserves the
manifold local linear structure while isometric mapping preserves
geodesic distance. Their inherent flexibility predisposes these methods
for extensions to functional data, where one rarely would have prior
information available about the nature of the underlying manifold.

Our goal is to explore manifold representations of functional data.
Which observed sets of functions are likely to lie on a low-dimensional
manifold? And how should this be taken into consideration? In contrast
to multivariate data, functional data are recorded on a time or
location domain, and commonly are assumed to consist of sets of smooth
random functions. Auspicious examples where functional manifold
approaches may lead to improved representations include time-warped
functional data \cite{wan99,ger04}, density functions~\cite{kne01}, and
functional data with pre-determined and interpretable modes~\cite{ize07}.
In such situations, the established linear functional
approaches, such as cross-sectional mean and functional principal
component analysis (FPCA) often fail to represent the functional data
in a parsimonious, efficient and interpretable way. Manifold approaches
are expected to be especially useful to represent functional data
inherently lying on a low-dimensional nonlinear space.

In this paper, we develop a framework for modeling $L^2$ functions on
unknown manifolds and propose pertinent notions, such as manifold mean,
manifold modes of functional variation and functional manifold
components, as elements of a functional manifold component analysis
(FMCA). Manifold means complement notions of a specifically modified
functional mean, such as the ``structural mean''~\cite{kne92}. A major
motivation for this proposal is that functional principal component
plots, for example, displaying second versus first component, are quite often
found to exhibit ``horseshoe'' shapes, that is, nonlinear dependence in
the presence of uncorrelatedness (as principal components by definition
are always uncorrelated). An example of this ``horseshoe shape'' is
provided by the Berkeley growth data (see upper right panel of Figure
\ref{fig5}). In such situations, one may wish to ``unwrap'' the
``horseshoe'' into linear structures by techniques similar to those
used in nonlinear dimension reduction. When attempting to ``unwrap''
functional data, one encounters specific difficulties: Often the
underlying smooth functions are not directly observed, but instead need
to be inferred from a limited number of noise-contaminated measurements
that contain the available information for each subject in the sample.
To address these problems, we develop a~modified ISOMAP \cite{ten00}
procedure, by adding a data-adaptive penalty to the empirical geodesic
distances, and employ local smoothing to recover the manifold.

The paper is organized in the following way. In Section \ref{sec2}, we
describe what we mean by a functional manifold, manifold mean, manifold
modes of functional variation and functional manifold components. We
develop corresponding estimates\vadjust{\goodbreak} in Section \ref{sec3} and discuss their
asymptotic properties in Section \ref{sec4}. Sections \ref{sec5} and
\ref{sec6} are devoted to illustrations of the proposed methodology for
both simulated and real data. Detailed proofs can be found in an online
supplement \cite{supp}.

%s2 #&#
\section{Manifolds in function space}\label{sec2}

%s2.1 #&#
\subsection{Preliminaries}\label{sec21}

A manifold $\mani$ can be expressed in terms of an atlas consisting of
a group of charts $(U_\alpha,\varphi_\alpha)$, where $U_\alpha$ are
open sets covering~$\mani$ and $\varphi_\alpha$, the coordinate maps,
map the corresponding $U_\alpha$ onto an open subset of $\R^d$.
Additional assumptions on $\varphi_\alpha$ are usually imposed in order
to study the structure of $\mani$ \cite{doc92,hel01}.

In this paper, we only consider ``simple'' functional manifolds $\mani
$ in $L^2$ space, where $\mani$ is \textit{isomorphic} to a subspace of
the Euclidean space, that is, the manifold can be represented by a
coordinate map $\varphi\dvtx \R^d\rightarrow\mani\subset L^2$, such
that $\varphi$ is bijective, and both $\varphi$, $\varphi^{-1}$ are
continuous, in the sense that if $\bs{\theta}_n,\bs{\theta}\in\R^d$ and
$\|\bs{\theta}_n-\bs{\theta}\|\rightarrow0$,
$\|\varphi(\bs{\theta}_n)-\varphi(\bs{\theta})\|_{L^2}\rightarrow0$; if
$x_n, x\in\mani$ and $\|x_n-x\|_{L^2}\rightarrow0$,
$\|\varphi^{-1}(x_n)-\varphi^{-1}(x)\|\rightarrow0$. Here, $d$ is the
intrinsic dimension of the manifold $\mani$. Such ``simple'' manifold
settings have been commonly considered in the dimension reduction
literature, for example in \cite{ten00}.

For a continuous curve defined on the manifold $\gamma
\dvtx[0,1]\rightarrow
\mathcal{M}$, define the length operator
%
%e2.1 #&#
%
\begin{equation}\label{211}
L(\gamma)={\sup\sum_{i=0}^{n-1}}\|\gamma(s_{i+1})-\gamma(s_i)\|_{L^2},
\end{equation}
where the supremum is taken over all partitions of the interval $[0,1]$
with arbitrary break points $0=s_0<s_1<\cdots<s_n=1$. We call $\varphi$
an isometric map if $L(\gamma)=L(\varphi^{-1}\circ\gamma)$ for any
continuous $\gamma$, where $L(\varphi^{-1}\circ\gamma)$ is similarly
defined as in (\ref{211}) with the $L^2$ norm replaced by the Euclidean
norm. We say~$\mani$ is an isometric manifold if there exists an
isometric coordinate map~$\varphi$. The isometry assumption is
pragmatically desirable and can be found in many approaches
\cite{ten00,don03}. Conditions under which isometry holds for image
data are discussed in \cite{don05}.

We use the notation $\psi\equiv\varphi^{-1}$ and refer to $\psi$ as the
representation map. The manifold $\mani$ is naturally equipped with the
$L^2$ distance, which, due to the nonlinearity of $\mani$, is not an
adequate metric \cite{ten00}. More useful is the geodesic distance
%
%e2.2 #&#
%
\begin{equation}\label{212}
d_g(x_1,x_2)=\inf\{L(\gamma)\dvtx\gamma(0)=x_1,\gamma(1)=x_2\},
\end{equation}
where the
infimum is taken over all continuous paths $\gamma$ on $\mani$. The geodesic
distance is the length of the shortest path on $\mani$ connecting the two
points, and therefore is adapted to $\mani$.

%s2.2 #&#
\subsection{Manifold mean and manifold modes of variation}\label{sec22}

Suppose $\mani$ is a~functional manifold of intrinsic dimension $d$ and
$\psi$ is a\vadjust{\goodbreak} representation map for $\mani$. Define, with respect to a
probability measure $Q$ in $\mathbb{R}^d$,
%
%e2.3 #&#
%
\begin{equation} \label{213}
\umani=\E\{\psi(X)\},\qquad \Xmu=\psi^{-1}(\umani),
\end{equation}
where $\umani$ is the mean in the $d$-dimensional representation
space, and
$\Xmu$ is the manifold mean in $L^2$ space. If $\mani$ is isometric, the
manifold mean $\Xmu$ is uniquely defined for all isometric representation
maps, as the following results shows.
%
%pr1 #&#
%
\begin{Prop}\label{prop1}
Suppose the random function $X$ lies on a functional manifold $\mathcal{M}$
of intrinsic dimension $d$ and $\psi$ is a representation map for~$\mani$. If
$\psi$ is isometric, the manifold mean $\Xmu$ in (\ref{213}) has the
following alternative expression:
%
%e2.4 #&#
%
\begin{equation} \label{215}
\Xmu= \argmin{x \in\mathcal{M}} \E d^2_g(x, X),
\end{equation}
where $d_g$ denotes the geodesic distance defined in (\ref{212}).
\end{Prop}

The expected value in equation (\ref{215}) is with respect to the
probability measure that is induced by the map $\varphi$; see also
\cite{bic07}. Equation (\ref{215}) defines the Fr\'{e}chet mean for
geodesic distance $d_g(\cdot,\cdot)$, and therefore does not depend
on the
choice of the isometric map $\psi$. The motivation to consider the manifold
mean is that the traditional cross-sectional mean for functional data in
$L^2$ has significant drawbacks as a measure of location when the data indeed
lie on a nonlinear functional manifold. Estimates of $L^2$ means,
obtained by
averaging observed sample curves, can be far away from the data cloud
in such
situations, and therefore do not represent the data in a meaningful way.
Going beyond the mean, one encounters analogous problems when linearly
representing such random functions in an $L^2$ basis, such as the
Fourier, B
spline or eigenfunction basis.

Consider random functions $X\in L^2(\mathcal{T})$ defined on a bounded
domain~$\mathcal{T}$. With $\mu(t)=\E X(t)$ and $G(t,s)=\mbox{Cov}(X(t),X(s))$,
according to Mercer's theorem \cite{ash75}, if the covariance function
$G(t,s)$ is jointly continuous in~$t$,~$s$, there is an orthonormal expansion
of $G(t,s)$ in terms of the eigenvalues $\{\lambda_k\dvtx\allowbreak k \geq1\}$ (ordered
nonincreasingly) and associated eigenfunctions \mbox{$\{\phi_k\dvtx k
\geq1\}$},
%
%e2.5 #&#
%
\begin{equation}\label{216} G(t,s)=\sum_{k=1}^\infty\lambda_k\phi
_k(t)\phi_k(s),\qquad
t,s\in\mathcal{T}.
\end{equation}
By the Hilbert--Schmidt theorem \cite{gre50,rie90},
$X$ can be expressed in terms of the so-called Karhunen--Lo\`{e}ve
representation,
%
%e2.6 #&#
%
\begin{eqnarray} \label{217}
X(t)&=&\mu(t)+\sum_{k=1}^\infty\xi_{k} \phi_k(t),\qquad
t\in\mathcal{T},\nonumber\\[-8pt]\\[-8pt]
\xi_k&=&\int_{\mathcal{T}}
\bigl(X(t)-\mu(t)\bigr)\phi_k(t)\,dt,\nonumber
\end{eqnarray}
where the $\xi_{k}$ are uncorrelated random variables with mean $0$ and
variance~$\lambda_k$, known as functional principal components (FPCs).

In the manifold case, the FPCs intrinsically lie on a $d$-dimensional
manifold. Therefore, we expect that the FPCs do not provide a parsimonious
representation of $X$. A better adapted and more compact representation can
be obtained through nonlinear manifold modes of functional variation
that are
defined below. The established eigenfunction-based modes of functional
variation \cite{cas86,jon92} are
%
%e2.7 #&#
%
\begin{equation} \label{219}
X_{j,\alpha}=\mu+\alpha\lambda_j^{1/2}\phi_j,\qquad j=1,2,
\ldots, \alpha\in\R,
\end{equation}
where factors $\lambda_j^{1/2}$ standardize the scale for different
$j$ and
the functional variation in the direction of eigenfunction $\phi_j$ is
visualized by the changing of functional shapes as $\alpha$ varies. However,
when the functional data lie on a~manifold, neither $\mu$ nor
$X_{j,\alpha}$
may belong to $\mani$, so that these linear modes will not provide a sensible
description of the variation in the data.

To address this problem, we define functional manifold component (FMC)
vectors $\ejmani\in\R^d$, $j=1,\ldots,d$, by the eigenvectors of the
covariance matrix of~$\psi(X)\in\R^d$, that is,
%
%e2.8 #&#
%
\begin{equation}\label{2110}
\Cov{(\psi(X))}=\sum_{j=1}^d\lambda_j^{\mathcal{M}}\ejmani
{\ejmani}^T,
\end{equation}
where $\lambda_1^{\mathcal{M}}\geq\cdots\geq\lambda_d^{\mathcal
{M}}$ are
the eigenvalues of $\Cov{(\psi(X))}$. The manifold modes of functional
variation are
%
%e2.9 #&#
%
\begin{equation} \label{2111}
X^{\mathcal{M}}_{j,\alpha} = \psi^{-1}\bigl(\umani+\alpha\slj
\ejmani\bigr),\qquad j=1,\ldots,d, \alpha\in\R,
\end{equation}
where $\umani$ is the mean in the $d$-dimensional representation space
according to measure $Q$, as given in (\ref{213}). A distinct
advantage of
manifold-based modes of functional variation over the principal component
based version (\ref{219}) is that in (\ref{2111}) only finitely many
modes are needed, while (\ref{219}) requires potentially infinitely many
components. The manifold modes $X^{\mathcal{M}}_{j,\alpha}$ are
unique for
the case of isometric~$\mani$, as shown in the following.\vspace*{-3pt}
%
%pr2 #&#
%
\begin{Prop}\label{prop2}
Suppose $\psi$ and $\tilde{\psi}$ are two isometric representation~maps
for a
functional manifold $\mathcal{M}$ of intrinsic dimension $d$. Let
$X^{\mathcal{M}}_{j,\alpha}$ be the~$j$th~ma\-nifold mode defined in
(\ref{2111}) based on representation map $\psi$,
and~$\tilde{X}^{\mathcal{M}}_{j,\alpha}$ be the $j$th manifold mode
using map
$\tilde{\psi}$. Then
$X^{\mathcal{M}}_{j,\alpha}=\tilde{X}^{\mathcal{M}}_{j,\alpha}$
for all
$\alpha\!\in\!\R$ and~\mbox{$1\leq j \leq d$}, if the eigenvalues of $\Cov
{(\psi(X))}$
and of $\Cov{(\tilde{\psi}(X))}$ are of multiplicity one.\vspace*{-3pt}
\end{Prop}

For each $X \in\mani$, given the representation map $\psi$, $X$ can be
uniquely represented (due to the bijectivity of $\psi$) by a vector
$\bs{\vartheta}=(\vartheta_1,\ldots,\vartheta_d)\in\R^d$
%
%e2.10 #&#
%
\begin{equation} \label{2112}
X = \psi^{-1}\Biggl(\umani+\sum_{j=1}^d \vartheta_j
\ejmani\Biggr),\qquad
\vartheta_j = \langle\psi(X)-\umani,\ejmani\rangle,\quad
j=1,\ldots,d,\hspace*{-42pt}\vadjust{\goodbreak}
\end{equation}
where $\umani$ and $\ejmani$ are defined in (\ref{213}) and
(\ref{2110}),
respectively, $\langle\cdot,\cdot\rangle$ is the inner product in
$\R
^d$ and
$\vartheta_j$ are uncorrelated r.v.s with mean $0$ and variance
$\lambda_j^{\mathcal{M}}$. We call $\vartheta_j$ the functional manifold
components (FMCs) in the representation space.

%s3 #&#
\section{Estimating functional manifolds}\label{sec3}

Suppose we observe $\{Y_{ij}\dvtx1\leq i \leq n;\allowbreak
1 \leq j \leq n_i\}$ which are noise-contaminated measurements made on
$n$ independent realizations $X_i$ of a random
function $X\in\mani$, according to the data model
\[
Y_{ij}=X_i(t_{ij})+\varepsilon_{ij}.
\]
Here the $t_{ij}$ are the time points where the functions are sampled, and
the~$\varepsilon_{ij}\in\R$ are i.i.d. errors with mean 0 and variance
$\sigma^2$. A first task is to find an approximation $\hat{\psi}$ to the
representation map $\psi$ based on the observed~$Y_{ij}$. We also
require the
inverse~$\hat{\psi}^{-1}$. Prior knowledge about the data may suggest a
specific form for $\psi$ \cite{ize07}, or one may have direct observations
of~$\psi(X_i)$. But in general, the representation map $\psi$ is
unknown and
needs to be determined from the data.

%s3.1 #&#
\subsection{Inferring $d$-dimensional manifold representations}\label{sec31}

Following \cite{ten00}, we use the pairwise distances between observed data
to obtain a map $\psi$ that preserves the geodesic distances. Alternative
approaches include LLE~\cite{row00} and Laplacian eigenmaps
\cite{bel03}. While these methods have been developed for multivariate data,
they can be adapted to functional data in a two-step procedure as follows.

In a first step, given an intrinsic dimension $d$ of $\mathcal{M}$,
adopt the
proposal of \cite{ten00} to obtain the function $\psi\dvtx L^2
\rightarrow
\R^d$ only at the sample points $\{X_1, \ldots, X_n\}$, where $X_i
\in L^2$,
by
%
%e3.1 #&#
%
\begin{equation} \label{31}
\hat{\psi}= \argmin{(\psi(X_1),\ldots,\psi(X_n))}
\sum_{i,j=1}^n\{\|\psi(X_i)-\psi(X_j)\|-d_g(X_i,X_j)\}^2.
\end{equation}
Here, $d_g(\cdot,\cdot)$ is the geodesic distance (\ref{212}) and the
minimum is taken over the vectors $\psi(X_i) \in\R^d, i=1,\ldots,n$,
formed by the values of $\psi$ on the functions~$X_i$, that is, the
goal is to find $n$ vectors $\hat{\psi}(X_i) \in\R^d, i=1,\ldots,n$,
that minimize~(\ref{31}). For this, one needs to estimate the geodesic
distances, and then the minimizer $\hat{\psi}(X_i)$ is obtained by
multidimensional scaling (MDS) based on estimates of
$d_g(X_i,X_j)$~\cite{cox01}. Our asymptotic results pertain to a~second
step, where the assumed smoothness of $\psi$ is invoked to obtain
global estimates for $\hat{\psi}$, as described in Section \ref{sec32}.
As for $\hat{\psi}(X_i), i=1,\ldots,n$, as determined by (\ref {31}),
we assume that the minimization in (\ref{31}) provides values on or
defines the target manifold at the sample points, that is, that
$\hat{\psi}(X_i)=\psi(X_i), i=1,\ldots,n$, or alternatively, that
$v_n=\hat{\psi}(X_i)-\psi(X_i)\rightarrow0$.

In order to approximate geodesic distances $d_g(X_i,X_j)$, we first aim at
estimates of the $L^2$ distances $\|X_i-X_j\|_{L^2}$. For this purpose, the
Karhunen--Lo\`{e}ve representation (\ref{217}) can be used to obtain fitted
curves,
%
%e3.2 #&#
%
\begin{equation}\label{312}
\hXiK(t)=\hat{\mu}(t)+\sumk\hat{\xi}_{ik}\hat{\phi}_k(t).
\end{equation}
Here,
$\hat{\mu}(t)$ and $\hat{G}(t,s)$ are first obtained by applying
local linear
one-dimen\-sional and two-dimensional smoothers to the pooled data; then
eigenfunctions $\hat{\phi}_k(t)$ and eigenvalues $\hat{\lambda}_k$ are
extracted by classical vector spectral analysis applied to a discretized
version of the estimate $\hat{G}(t,s)$ of the covariance surface
${G}(t,s)=\Cov(X(t),X(s))$; and then the FPCs $\xi_{ik}$ are
approximated by
discretizing integrals
%
%e3.3 #&#
%
\begin{equation}\label{314}
\hat{\xi}_{ik}=\sumj\{Y_{ij}-\hutij\}\hphik(t_{ij})(t_{ij}-t_{i,j-1})
\end{equation}
or
alternatively by conditional expectation (for details on these steps,
see~\cite{yao05}),
%
%e3.4 #&#
%
\begin{equation}\label{315}
\hat{\xi}_{ik}=\hat{\lambda}_k\hat{\bs{\phi}}^T_{ik}\hat{\Sigma
}^{-1}_{\mf{Y}_i}(\mf{Y}_i-\hat{\bs{\mu}}_i),
\end{equation}
where\vspace*{1pt} $\hat{\bs{\phi}}_{ik}=(\hat{\phi}_k(t_{i1}),
\ldots,\hat{\phi}_k(t_{in_i}))$,
$(\hat{\Sigma}_{\mf{Y}_i})_{jl}=\hat{G}(t_{ij},t_{il})+\hat{\sigma
}^2\mf
{I}$,
$1\leq j$, $l \leq n_i$,
$\hat{\bs{\mu}}_i=(\hat{\mu}(t_{i1}),\ldots,\hat{\mu
}(t_{in_i}))$, and
$\sigma^2$ is estimated from the difference between empirical
variances of
$Y_{ij}$ and $\hat{G}(t,s)$. The conditioning method (\ref{315}) is the
only~available option if the data are sparsely sampled. To ensure that
a~sufficiently large number of components is included in the truncated
expansion~(\ref{312}), one may choose $K$ by requiring a large
fraction of
variance explained (FVE), that is,
%
%e3.5 #&#
%
\begin{equation}\label{313}
K=\min_k\biggl\{k\dvtx\frac{\sum_{l\leq k}\hat{\lambda}_l}{\sum
_{l=1}^\infty
\hat{\lambda}_l}\geq1-\alpha\biggr\}
\end{equation}
for, say, $\alpha=0.05$, where the
$\hat{\lambda}_l$ are estimates of the eigenvalues $\lambda_l$ in
(\ref{216}). The resulting $L^2$ distances are
$\|\hXiK-\hat{X}_j^K\|_{L^2}=\{\sum_{k=1}^K(\hat{\xi}_{ik}-\hat
{\xi
}_{jk})^2\}^{1/2}$.

Note that alternatively to representation (\ref{312}), one can also
directly apply local constant or local linear smoothing to obtain smooth
trajectories in the case of dense and balanced designs, for example,
using Nadaraya--Watson kernel estimators,
%
%e3.6 #&#
%
\begin{equation}\label{311}
\tXit=\frac{\sum_{j=1}^{n_i}\kappa_1(h_1^{-1}(t_{ij}-t))Y_{ij}}
{\sum_{j=1}^{n_i}\kappa_1(h_1^{-1}(t_{ij}-t))},
\end{equation}
where $\kappa_1$
and $h_1$ are smoothing kernel and bandwidth. For the smoothing kernel one
can use any standard kernel such as the standard Gaussian density
function or
the Epanechnikov kernel, while in practice $h_1$ may be chosen by
cross-validation or generalized cross-validation.

Then the pairwise $L^2$ distances are simply
$\|\tilde{X}_i-\tilde{X}_j\|_{L^2}$. We will not explicitly explore this
alternative smoothing approach in our theoretical analysis, but note that
essentially the same results as those reported below hold for this
alternative approach, by minor extensions of our arguments. In the
implementations (simulation and data analysis), we use both
approaches~(\ref{312}) and (\ref{311}). The estimated random trajectories, obtained
though (\ref{312}) or (\ref{311}), generally are not lying on the
manifold $\mani$, as they are merely approximations to the true unknown
functions, due to additional noise and discrete sampling of the random
trajectories. However, these estimates, owing to their consistency, will
fall inside a small $L^2$-neighborhood around $\mani$. Asymptotic properties
are discussed in Section \ref{sec4}.

Since the geodesic is the shortest path connecting points on a
manifold, one
may first connect the points inside small $L^2$ neighborhoods and then define
the path between two far away points by moving along these small
neighborhoods, and then find the geodesic by the shortest path connecting
through such neighborhoods. This is essentially the idea of the ISOMAP
algorithm \cite{ten00}. The performance of this method however proved
somewhat unstable in our applications, as functional data typically
must be
inferred from discretized and noisy observations of underlying smooth
trajectories and therefore do not exactly lie on the manifold, as is assumed
in ISOMAP.

In such situations, due to random scattering of the data around the manifold,
the shortest path found by the ISOMAP criterion may pass through ``empty
areas'' outside the proper data cloud. This problem can be effectively
addressed by modifying the ISOMAP criterion, by additionally penalizing
against paths that include sections situated within ``empty regions'' with
few neighboring data points. Density-penalized geodesics are
characterized by
sequences of $L^2$ functions $(W_1, W_2, \ldots, W_m)$ from the starting
point $W_s=W_1$ to the end point $W_e=W_m$ of the geodesic, where each
of the
$W_j$ stands for one of the observed functions $X_i$ (with unrelated index),
and are defined as
%
%e3.7 #&#
%
\begin{eqnarray}\label{def43}
S(W_s,W_e)&=&\argmin{W_2,\ldots,W_{m-1}}\Biggl\{\sum
_{i=1}^{m-1}\|W_i-W_{i+1}\|_{L^2} \bigl(1+P_\delta(W_i,W_{i+1})\bigr)\dvtx
\nonumber\\[-8pt]\\[-8pt]
&&\hspace*{144pt}\|W_i-W_{i+1}\|_{L^2} <\varepsilon\Biggr\}. \nonumber
\end{eqnarray}
Here the parameter $\varepsilon$ limits the step length, and the penalty
function $P_\delta$ is determined by the density of the data cloud around
$W_i$ and $W_{i+1}$,
\[
P_\delta(W_i, W_{i+1})=\rho_{i,i+1}^{-2}I(\rho_{i,i+1}<\delta),
\]
where $\rho_{i,i+1}=\min\{\#\{W_j\dvtx\|W_j-W_i\|_{L^2}<\varepsilon\},
\#\{W_j\dvtx\|W_j-W_{i+1}\|_{L^2}<\varepsilon\}\}$ and $\#$ denotes the
cardinality
of a set. By selecting the parameter $\delta$, one can control the threshold
of the local density of points, below which the penalty~$P_\delta$
kicks in.
The ISOMAP algorithm corresponds to the special case where $\delta=0,
P_{\delta}=0$.

The choice $\delta>0$ leads to ``penalized ISOMAP'' or P-ISOMAP, where the
penalty parameter $\delta$ may be selected data-adaptively by
cross-validation. The choice of $\delta$ and also of the step size parameter
$\varepsilon$ is discussed in Section \ref{sec33}. If the manifold is very
smooth, a
large $\varepsilon$ and small $m$ will lead to a sufficiently good
estimate of
the geodesic distance. A detailed discussion of the convergence of the
estimated geodesics in the framework of ISOMAP can be found at
\url{http://isomap.stanford.edu/BdSLT.pdf}. For the proposed P-ISOMAP,
we implement
the minimization of $S(W_s,W_e)$ by Dijkstra's algorithm, which selects $m$
and the geodesic paths $(W_s=W_1,W_2,\ldots,W_{m-1},W_e=W_m)$. The resulting
estimated geodesic distance is
%
%e3.8 #&#
%
\begin{equation}\label{geo} \hat{d}_g(W_s,W_e)
=\sum_{j=1}^{m-1} \|\hat{W}_j-\hat{W}_{j+1}\|_{L^2},
\end{equation}
where
$\hat{W}_j=\tilde{W}_j$ or $\hat{W}^K_j$, depending on which preliminary
approximation is used for $W_j$.
Once these distances have been determined, an application of MDS yields
$\hat{\psi}(X_i)$, in the same way as in the standard ISOMAP method.

%s3.2 #&#
\subsection{Obtaining the global map and representing sample
trajectories}\label{sec32}

For any location $\bs{\theta} \in\R^d$, we find
$\hat{\psi}^{-1}(\bs{\theta})$ by local weighted averaging, that is,
%
%e3.9 #&#
%
\begin{equation}\label{332} \hat{\psi}^{-1}(\bs{\theta})=\frac
{\sumii
\kappa(H^{-1}(\hat{\psi}(X_i)-\bs{\theta})) \hat
{X}_i}{\sumii
\kappa(H^{-1}(\hat{\psi}(X_i)-\bs{\theta}))},
\end{equation}
where $\kappa$ is a $d$-dimensional kernel, like the Epanechnikov kernel
$\kappa(u_1,\ldots,\allowbreak u_d)=(\frac{3}{4})^d\prod_{k=1}^d\{(1-u_k^2)\mf
{I}(|u_k|<1)\}$, with $H=h\mf{I}_{d\times d}$ for a suitably chosen
bandwidth $h$, $\hat{X}_i$
could be either $\tilde{X}_i$ as in (\ref{311}) or $\hat{X}^K_i$
as in
(\ref{312}), and $\hat{\psi}(X_i)$ is defined after (\ref{geo}).
We use
cross-validation to select $h$ (see Section \ref{sec33}). The
asymptotic properties
of (\ref{332}) will be discussed in Section \ref{sec4}.

Specifically, as predictor of $X_i$, we propose
%
%e3.10 #&#
%
\begin{equation}\label{pred}
\hat{X}_i^\mani=\frac{\sum_{j\neq i}
\kappa(H^{-1}(\hat{\psi}(X_i)-\hat{\psi}(X_j)))
\hat{X}_j}{\sum_{j\neq i}
\kappa(H^{-1}(\hat{\psi}(X_i)-\hat{\psi}(X_j)))},
\end{equation}
borrowing
strength from local neighbors in the $d$-dimensional representation space.
This can be seen as an alternative to representation (\ref{312}),
where we
use the FPCs and borrow strength from the whole data set to estimate
functional mean and eigenbasis. As before, we note that (\ref{pred})
is not
necessarily in $\mani$,\vadjust{\goodbreak} but will be in a small neighborhood asymptotically
and in comparison with (\ref{312}), (\ref{pred}) usually proves to
be a
much better predictor of $X_i$ for functional manifold data as shown in the
simulations and applications in Section \ref{sec5}. Asymptotic
properties are
discussed in Section \ref{sec4}.

Definition (\ref{213}) suggests to estimate the manifold mean by
%
%e3.11 #&#
%
\begin{equation}\label{341} \hXmu=\frac{\sumii
\kappa(H^{-1}(\hat{\psi}(X_i)-\humani)) \hat
{X}_i}{\sumii
\kappa(H^{-1}(\hat{\psi}(X_i)-\humani))},
\end{equation}
where
$\humani=\frac{1}{n}\sumii\hat{\psi}(X_i)$. The FMC vectors
$\ejmani$ defined
in (\ref{2110}) are estimated by eigendecomposition of the sample
covariance matrix of $\hat{\psi}(X_i)$, that is,~$\hat{\lambda
}^{\mathcal{M}}_j$
and $\hejmani$ are such that
%
%e3.12 #&#
%
\begin{eqnarray} \label{342}
\sum_{j=1}^d \hat{\lambda}_j^{\mathcal{M}} \hejmani\hejmani^T&=&
\frac
{1}{n-1}\Biggl\{\sum_{i=1}^n \hat{\psi}(X_i)
\hat{\psi}^T(X_i)\nonumber\\[-8pt]\\[-8pt]
&&\hspace*{30.8pt}{}-\frac{1}{n}\Biggl(\sum_{j=1}^n
\hat{\psi}(X_j)\Biggr)\Biggl(\sum_{j=1}^n \hat{\psi}(X_j)
\Biggr)^T\Biggr\},
\nonumber
\end{eqnarray}
where\vspace*{2pt} the $\hat{\lambda}_j^{\mathcal{M}}$ are ordered to be
nonincreasing in
$j$. From (\ref{2111}) and (\ref{332}), we obtain estimates of the
manifold modes as
%
%e3.13 #&#
%
\begin{equation}\label{343} \hXj=\frac{\sumii
\kappa(H^{-1}\{\hat{\psi}(X_i)-\humani-\alpha
(\hat{\lambda}_j^{\mathcal{M}})^{1/2} \hejmani\}) \hat
{X}_i}{\sumii
\kappa(H^{-1}\{\hat{\psi}(X_i)-\humani-
\alpha(\hat{\lambda}_j^{\mathcal{M}})^{1/2} \hejmani\})},
\end{equation}
where $j=1,\ldots,d$ and $\alpha\in\R$.

%s3.3 #&#
\subsection{Selection of auxiliary parameters}\label{sec33}

We use $10$-fold cross-validation to simultaneously choose the step size
$\varepsilon$, the truncation parameter $\delta$, and the smoothing bandwidth
$h$ (see Sections \ref{sec31} and \ref{sec32}). The number of
candidates for
$\varepsilon$
and $\delta$ is kept small so that the cross-validation procedure runs
reasonably fast. Candidates for the step size $\varepsilon$ are the median
distance of the $5$th, the $8$th and the $12$th nearest neighbor; those for
$\delta$ are selected such that $0\%$, $2\%$, $5\%$ and $10\%$ of the data
with the lowest local density estimates are penalized. Each of 10 subgroups
of curves denoted by $V_1,\ldots,V_{10}$ is used as a validation set,
one at
a time, while the remaining data are used as training set.

In an initial step, we use the whole data set and a given $\varepsilon$,
$\delta$ to determine~$\hat{\psi}(X_i)$, followed by estimation of
$X_i=\psi^{-1}(\bs{\vartheta}_i)$ for $X_i$ in the validation set, using
(\ref{332}) and assuming that only those $\hat{X}_j$ in the
training set
are known. Denoting the value of the estimated trajectory $X_i$,
evaluated at
time $t_{il}$, by~$\hat{X}_{il}$, the sum of squared prediction errors for
the validation set~$V_k$ is $\mathrm{SSPE}_k=\sum_{i\in
\mathrm{V}_k}\sum_{l=1}^{n_i}(\hat{X}_{il}-Y_{il})^2$, where
$Y_{il}=X_i(t_{il})+\varepsilon_{il}$ is the observed value of
trajectory $X_i$
at time $t_{ij}$. The cross-validation choice is the minimizer of
$\mathrm{MSPE}(\varepsilon,h,\delta)=\frac{\sum_{k=1}^{10}
\operatorname{SSPE}_k}{\sum_{i=1}^n n_i}$.\vadjust{\goodbreak}

Following \cite{ten00}, the intrinsic dimension $d$ can be chosen by the
$1-\beta$ fraction of distances explained (FDE), that is,
%
%e3.14 #&#
%
\begin{equation}\label{331}
d=\min_p\biggl\{p\dvtx\frac{\|\hat{D}^p-D\|_F}{\|D\|_F}<\beta\biggr\},
\end{equation}
where we choose
$\beta=0.05$ and $D$, $\hat{D}^p$ are $n$ by $n$ distance matrixes with
$D_{ij}=\hat{d}_g(X_i,X_j)$ as in (\ref{geo}),
$\hat{D}^p_{ij}=\|\hat{\psi}^p(X_i)-\hat{\psi}^p(X_j)\|$\vspace*{1pt} and where
$\hat{\psi}^p$ denotes the MDS solution (\ref{31}) in $\R^p$, and
\mbox{$\|\cdot\|_F$} is the matrix Frobenius norm,
$\|D\|_F=\{\sum_{i,j}D^2_{ij}\}^{1/2}$. Note that $\|\hat{D}^p-D\|_F$
is the square root of the minimized value of (\ref{31}).

%s4 #&#
\section{Asymptotic properties}\label{sec4}

We provide the specific convergence rate\break of~$\hat{X}_i^K$, defined in
(\ref{312}), under assumptions (A1)--(A5) in the \hyperref
[app]{Appendix}. Note that
condition (A3) requires that the random functions are sampled at a~dense
design. Our starting point is that the manifold can be well identified
at the
sample points through ISOMAP, or alternatively, that the ISOMAP identified
manifold may be viewed as the target. The difference between the target and
the identified manifold from ISOMAP is quantified by a rate $v_n$ that is
assumed as given; if the target manifold corresponds to the manifold as
identified at the sample points, we may set $v_n=0$. The theoretical analysis
aims to justify the new manifold representations that we propose, and for
this it is essential to consider the behavior of the estimates across the
entire function space. Therefore, our theoretical results demonstrate
how to
extend local behavior at the sample points to obtain global consistency of
the proposed functional manifold representations.

As the convergence is for $K=K_n\rightarrow\infty$ as $n\rightarrow
\infty$,
the rate of decline of the eigenvalues in (\ref{216}) and also lower bounds
on the spacing of consecutive eigenvalues, as postulated in (A4) are
relevant, with a requirement of polynomially fast declining eigenvalues.
Required smoothness and boundedness assumptions for $X\in\mani$ are
as in
(A5).
%
%pr3 #&#
%
\begin{Prop} \label{prop3}
Assume \textup{(A1)--(A5)} in the \hyperref[app]{Appendix}, and define
$r_n=\max\{\frac{1}{\sqrt{n}h^2_G},\frac{1}{\sqrt{n}h_{\mu}},
\frac{1}{\sqrt{n}h_V}\}$. If there are infinitely many nonzero eigenvalues
$\lambda_k$ in (\ref{216}), which are all of multiplicity one, then for
sequences $K=K_n\rightarrow\infty$, subject to $r_n
K^{\alpha_2+{1/2}}\rightarrow0$, where $\alpha_2$ is a
constant such
that $\lambda_{k}-\lambda_{k+1}>C_2 k^{-\alpha_2}$ for some $C_2>0$
and where
$K\leq K_0$ with $K_0=\min\{i\dvtx\lambda_i-\lambda_{i+1}\leq2D_n\}-1$ and
$D_n=\{\int_{\mathcal{T}^2} (\hat{G}(t,s)-G(t,s))^2\,dt\,ds\}^{1/2}$ where
$G$ is
defined in~(\ref{216}) and $\hat{G}$ is defined after (\ref
{312}), it
holds that
%
%e4.1 #&#
%
\begin{equation}\label{prop11}
\|\hat{X}_i^K-X_i\|_{L^2}=O_p\bigl(r_n K^{\alpha_2+{1/2}} +
K^{-(\alpha_1-1)/2}\bigr)
\end{equation}
for $\hat{X}_i^K$ defined in (\ref{312}), where $\alpha_1$ is such that
$\lambda_k<C_1 k^{-\alpha_1}$ for all $k$ and some $C_1<\infty$.
\end{Prop}

We note that under the assumptions, $K_0 \rightarrow\infty$. The
first term
on the r.h.s. of (\ref{prop11}) is due to estimation error and the second
term is due to truncation error. In the special case when there are only
finitely many nonzero $\lambda_k$ in~(\ref{216}), it can be shown
that the
rate in (\ref{prop11}) simply becomes $O_p(r_n)$. Next we discuss the
convergence of the estimates that appear in (\ref{342}).\vspace*{-3pt}
%
%pr4 #&#
%
\begin{Prop} \label{prop4}
Under \textup{(B1)} and \textup{(B2)} in the \hyperref[app]{Appendix},
%
%e4.2 #&#
%
\begin{equation}\label{prop41}
\|\humani-\umani\|=O_p\biggl(v_n+\frac{1}{\sqrt{n}}\biggr),
\end{equation}
where $\umani$ and $\humani$ are defined in (\ref{213}) and (\ref{341}),
and $v_n=\sup_{i=1,\ldots,n}\|\hat{\psi}(X_i)-\psi(X_i)\|$. If the $j$th
eigenvalue of $\Cov(\psi(X))$ is of multiplicity one, then
%
%e4.3 #&#
%
\begin{eqnarray}\label{prop42}
\|\hejmani-\ejmani\|&=&O_p\biggl(v_n+\frac{1}{\sqrt{n}}\biggr),
\\[-3pt]
%
%e4.4 #&#
%
\label{prop43}
|\hat{\lambda}^\mani_j-\lambda^\mani_j|&=&O_p\biggl(v_n+\frac{1}{\sqrt{n}}\biggr),
\end{eqnarray}
where $\lambda^\mani_j$, $\ejmani$, $\hejmani$ and $\hat{\lambda
}^\mani_j$
are defined in (\ref{2110}) and (\ref{342}), respectively.\vspace*{-3pt}
\end{Prop}
%
%th1 #&#
%
\begin{Thm} \label{thm1}
Under \textup{(A1)--(A5)}, \textup{(B1)}, \textup{(B2)} and \textup{(C1)--(C3)}
in the \hyperref[app]{Appen-} \hyperref[app]{dix}, assume that the
density function $f$ of $\psi(X)\in\R^d$ satisfies $f(\bs{\theta})>0$
for a~specific $\bs{\theta}=\psi(x)$ and that $h>0$ is selected such that
$h\rightarrow0$, $n^{-1}h^{-2(d+1)}\rightarrow0$ and $h^{-(d+1)}\E
v_n\rightarrow0$. Then $\hat{\psi}^{-1}(\bs{\theta})$ defined in
(\ref{332}), using $\hat{X}_i=\hat{X}_i^K$, is a consistent
estimate of
$\psi^{-1}(\bs{\theta})$. Specifically, defining
$T^K_{\bs{\phi}}=\{\sum_{k>K}\xi^2_{k}\}^{1/2}$ where $\xi
_k=\int
(X-\E
X)\phi_k$ and the orthonormal basis $\{\phi_k\dvtx k\geq1\}$ is given
in~(\ref{216}), and defining
$R_K(\bs{\theta})=T^K_{\bs{\phi}}(\psi^{-1}(\bs{\theta}))$, where
$R_K(\bs{\theta})\rightarrow0$ as $K=K_n\rightarrow\infty$, it
holds that
%
%e4.5 #&#
%
\begin{eqnarray}\label{thm11}
&&
\|\hat{\psi}^{-1}(\bs{\theta})-\psi^{-1}(\bs{\theta})\|_{L^2}\nonumber\\[-8pt]\\[-8pt]
&&\qquad=O_p\biggl(h^2+\frac{1}{\sqrt{nh^d}}+\frac{v_n}{h}+R_K(\bs{\theta
})+K^{\alpha
_2+{1/2}}r_n\biggr),\nonumber
\end{eqnarray}
where $r_n$, $\alpha_2$ and $v_n$ are as in assumptions \textup{(A3)},
\textup{(A4)} and \textup{(B1)}.\vspace*{-3pt}
\end{Thm}

Note that $R_K(\bs{\theta})$ corresponds to the truncation error for
$\psi^{-1}(\bs{\theta})\in\mani$. The last term $K^{\alpha
_2+{1/2}}r_n$
is due to the estimation error as in Lemma 1. The middle term
$\frac{v_n}{h}$ reflects the estimation error of the weights, which is
influenced by the scale of the bandwidth. The first part
$h^2+\frac{1}{\sqrt{nh^d}}$ is the optimal rate when the $X_i$ and
$\psi
$ are
known, reflecting an intrinsically $d$-dimensional smoothing problem. Related
findings are discussed in \cite{bic07}.

For the manifold modes, we obtain the following corollary.\vspace*{-3pt}
%
%co1 #&#
%
\begin{Cor} \label{Cor1}
Under the conditions of Theorem \ref{thm1}, for a given $\alpha\in\R$ and
$1\leq
j\leq d$, assume that $f(\umani+\alpha\slj\ejmani)>0$ and that $h$
is chosen
as in Theorem \ref{thm1}. Then the estimated manifold\vadjust{\goodbreak} modes $\hXj$ as in
(\ref{343}), substituting $\hat{X}_i=\hat{X}^K_i$, are consistent.
Specifically,
%
%e4.6 #&#
%
\begin{eqnarray}\label{cor11}
&&
\|\hXj-\Xj\|_{L^2}\nonumber\\[-8pt]\\[-8pt]
&&\qquad=O_p\biggl(h^2+\frac{1}{\sqrt{nh^d}}+\frac
{v_n}{h}+\frac
{1}{\sqrt{n}h}+R_K+K^{\alpha_2+{1/2}}r_n\biggr),\nonumber
\end{eqnarray}
where $R_K=T^K_{\bs{\phi}}(\Xj)$.
\end{Cor}

An immediate consequence of these results is that the manifold representation
given in (\ref{2112}) provides a consistent representation of all random
functions in the functional manifold. Proofs of all propositions,
theorem and corollary can be found in the supplementary file \cite{supp}.

%s5 #&#
\section{Examples and simulation study}\label{sec5}

%s5.1 #&#
\subsection{Functional manifolds and isometry}\label{sec51}

To illustrate our methods and to discuss the impact of the critical isometry
assumption, we consider the following three example functional manifolds:
\begin{longlist}[(iii)]
\item[(i)] A one-dimensional $(d=1)$ functional manifold
\begin{eqnarray*}
\mani_1&=&\biggl\{X\in L^2([-4,4])\dvtx
X(t)=\mu(h_{\alpha}(t)), \\
&&\hspace*{7pt} h_{\alpha}(t)=\frac{8\int_0^{t/8+0.5} s^{\alpha}(1-s)\,ds}{\int_0^1
s^{\alpha}(1-s)\,ds}-4, \alpha>-1\biggr\},
\end{eqnarray*}
where $\mu(t)=\frac{2}{\sqrt{\pi}}\exp\{-\frac{1}{2}(t+2)^2\}
+\frac
{1}{\sqrt{2\pi}}\exp\{-2(t-2)^2\}$.
This corresponds to random warping of a common shape function
$\mu$, which has two peaks. The time warping function $h_{\alpha}$
is generated from the cumulative Beta distribution family and
$\alpha$ is a random parameter, $\alpha=\max(-1,Z)$, where $Z \sim
\mathrm{N}(0,0.09)$.
\item[(ii)] A two-dimensional $(d=2)$ functional manifold
\begin{eqnarray*}
\mani_2&=&\biggl\{X\in L^2([-4,4])\dvtx
X(t)=\frac{1}{\sqrt{2\pi\alpha^2}}\exp\biggl[-\frac{1}{2\alpha^2}
(t-\beta)^2\biggr],\\
&&\hspace*{194pt}\alpha>0, \beta\in\R\biggr\}.
\end{eqnarray*}
This manifold is a
collection of Gaussian densities, corresponding to a shift-scale
family, where $\alpha=\max(0,Z)$, $Z \sim\mathrm{N}(1,0.04)$ and
$\beta\sim\mathrm{N}(0,1)$.
%
%f1 #&#
%
\begin{figure}

\includegraphics{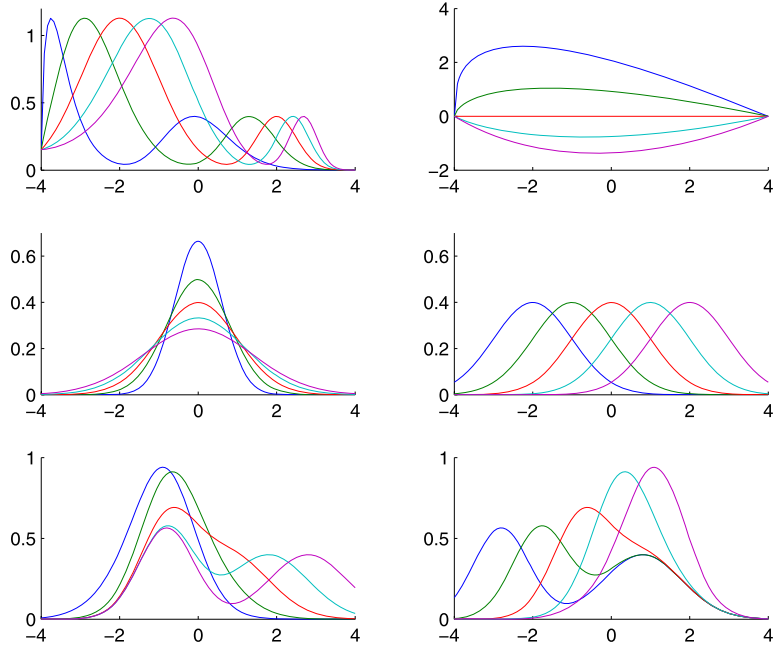}

\caption{Manifolds $\mani_1$--$\mani_3$. Top left panel:
functions on $\mani_1$ for $\alpha=0.6, 0.8, 1.0, 1.2, 1.4$. Top
right panel: corresponding identity-subtracted warping functions
$h_{\alpha}(t)-t$. Middle left panel: functions on $\mani_2$ for
$\alpha=0.4, 0.7, 1.0, 1.3, 1.6$ and $\beta=0$. Middle right panel:
functions on $\mani_2$ for $\beta=0.4, 0.7, 1.0, 1.3, 1.6$ and
$\alpha=0$. Bottom left panel: functions on $\mani_3$ for
$\alpha=-2, -1, 0, 1, 2$ and $\beta=0$. Bottom right panel:
functions on $\mani_3$ for $\beta=-2, -1, 0, 1, 2$ and $\alpha=0$.}
\label{fig1}
\vspace*{3pt}
\end{figure}
\item[(iii)] Another two-dimensional $(d=2)$ functional manifold
\begin{eqnarray*}
&&\mani_3=\biggl\{X\in L^2([-4,4])\dvtx X(t) =\frac{1}{\sqrt{2\pi
}}\exp\biggl\{
-\frac{1}{2}(t-0.8-\alpha)^2\biggr\}\\
&&\hspace*{37.5pt}\hspace*{112pt}{}
+\frac{1}{\sqrt{\pi}}\exp\{-(t+0.8-\beta)^2\},
\alpha,\beta
\in
\R\biggr\},
\end{eqnarray*}
a mixture of two peaks with randomly varying centers, where $\alpha
\sim
\mathrm{N}(0,1)$ and $\beta\sim\mathrm{N}(0,1)$. Note that the two peaks
will merge to a larger peak when their locations are close, so this set
of functions has a randomly varying number of peaks.
\end{longlist}

Functional manifolds $\mani_1$--$\mani_3$ are illustrated in Figure \ref
{fig1}. We note
that~$\mani_1$ is an isometric manifold and $\mani_2$ is approximately
isometric, while $\mani_3$ is not isometric. This can be seen as
follows. For
functions $X \in L^2$ on a differentiable isometric manifold with
representation $X=\psi^{-1}(\theta_1,\ldots,\theta_d)$, using the definition
of isometry given after (\ref{211}), the condition
${\int_{\theta^0_k}^{\theta^1_k}}\|\frac{\partial X}{\partial
\theta_k}(t)\|_{_{L^2}}\,d\theta_k\equiv\theta^1_k-\theta^0_k$ for
$k=1,\ldots,d$
and any $\theta^0_k,\theta^1_k\in\R$ is equivalent to isometry. Therefore,
the existence of a parametrization of the map $\psi$ for which the $L^2$
norms of the partial derivatives of $X$ with respect to the parameter
components are constant is sufficient and necessary for $\psi$ to be
isometric. For one-dimensional manifolds such as $\mani_1$, one can always
find such a parametrization, as long as~$X$ is differentiable in the
parameter and the derivative is $L^2$ integrable in~$t$.\looseness=-1

For $\mani_2$, such a parametrization does not exist, but since
$\|\frac{\partial X}{\partial\alpha}(t)\|_{L^2}=\frac{1}{\alpha
}c_1$ and
$\|\frac{\partial X}{\partial\beta}(t)\|_{L^2}=\frac{1}{\alpha
}c_2$ for
constants $c_1$, $c_2$ and as $\alpha$ is chosen to remain very close
to~$1$,
the natural parametrization approximately satisfies the condition for
isometry. In contrast to $\mani_1$ and $\mani_2$, the functional manifold
$\mani_3$ is nonisometric and we include it as an example how the proposed
methodology is faring when the key assumption of isometry is violated. As
our considerations take place in a manifold learning framework, where the
underlying manifold is unknown, an interesting aspect is to devise a
data-based check to gauge the degree to which the isometry assumption
can be
expected to be satisfied. A natural metric for such a check is the fraction
of distances explained (FDE), defined in (\ref{331}). This criterion
quantifies the percentage of geodesic distance that is preserved when fitting
a $d$-dimensional isometric manifold to the data. For cases where the
underlying manifold is actually nonisomorphic, the fitted manifold is an
isometric approximation to the true underlying manifold, obtained by
minimizing the stress function in the MDS algorithm.

An informal goodness-of-fit criterion for isometry is to require FDE to be
larger than 95\%, and choosing the manifold with the smallest dimension that
satisfies this criterion. In Table \ref{table1}, values for FDE
obtained for the
simulated data for manifolds $\mani_1$--$\mani_3$ under two signal-to-noise
ratios $R$ (defined in the following subsection) are reported, with dimension
$d$ ranging from $1$ to $5$. The well-known fact that the stress function
declines when the dimension of the projection space is increased underlies
the traditional MDS-Scree Plot \cite{cox01} and is reflected by the
observed increase in the values for FDE as dimension increases.

%t1 #&#
%
\begin{table}
\caption{Fraction of distances explained (\protect\ref{331}) for isometric
manifold fits with different dimension $d$ (other~parameters are optimized),
for two signal-to-noise ratios $R$}\label{table1}
\begin{tabular*}{\tablewidth}{@{\extracolsep{\fill}}lcd{1.4}cccc@{}}
\hline
& & \multicolumn{5}{c@{}}{$\bolds{d}$} \\[-4pt]
& & \multicolumn{5}{c@{}}{\hrulefill} \\
\textbf{Manifold} & \multicolumn{1}{c}{$\bolds{R}$}
& \multicolumn{1}{c}{\textbf{1}} & \multicolumn{1}{c}{\textbf{2}}
& \multicolumn{1}{c}{\textbf{3}} & \multicolumn{1}{c}{\textbf{4}}
& \multicolumn{1}{c@{}}{\textbf{5}} \\
\hline
$\mani_1$ & 0.1 & 0.998 & 0.999 & 0.999 & 0.999 &
0.999\\
& 0.5 & 0.9778 & 0.993 & 0.995 & 0.995 & 0.996 \\
[3pt]
$\mani_2$& 0.1 & 0.914 & 0.988 & 0.994 & 0.996 &
0.996 \\
& 0.5 & 0.902 & 0.971 & 0.974 & 0.978 & 0.980 \\
[3pt]
$\mani_3$ & 0.1 & 0.699 & 0.932 & 0.957 & 0.977 &
0.980 \\
& 0.5 & 0.639 & 0.906 & 0.948 & 0.955 & 0.958 \\
[3pt]
Growth & & 0.947 & 0.972 & 0.980 & 0.985 & 0.988 \\
Yeast & & 0.891 & 0.949 & 0.981 & 0.983 & 0.984 \\
Mortality & & 0.878 & 0.954 & 0.973 & 0.980 & 0.982 \\
\hline
\end{tabular*}
\end{table}

Applying the above check for isometry, we find that indeed the
dimensions of
the isometric manifold $\mani_1$ and the near-isometric manifold
$\mani_2$
are correctly selected, while the first two dimensions of the isometric
manifold approximation to the nonisometric manifold $\mani_3$ are not
sufficient. Thus, the nonisometric nature of $\mani_3$ means that the
dimension of the underlying functional manifold cannot be correctly
identified and instead the proposed algorithm will find a higher-dimensional
isometric manifold to represent $\mani_3$. The price to pay for a suitable
isometric approximation is increased dimensionality, which in this example
ends up larger than~2 for the approximating isometric manifold. We note that
an approximating isometric manifold can always be found, since the
linear and
therefore intrinsically isometric manifold of infinite dimensionality
that is
spanned by the eigenfunction basis contains the random functions of the
sample, according to the Karhunen--Lo\`eve theorem, and is always applicable.

While we can always find a near-isometric manifold of large enough
dimensionality with the proposed algorithm, when the data lie on a
lower-dimensional nonisometric manifold, these approximating isometric
manifolds may not be efficient, since they do not provide the
lowest-dimensional possible description of the data. Nevertheless, an
approximating isometric nonlinear manifold obtained by the proposed approach
often will present a much improved and lower-dimensional description when
compared to the alternative of classical linear basis representation.
This is
exemplified by the functional nonisometric manifold $\mani_3$, which
in the
following subsection is shown to be much better represented by an isometric
manifold than by a linear basis. So the price that the isometry assumption
exacts in nonisometric situations is that the proposed approach leads
to a
more or less suboptimal representation, which however will often be
substantially lower-dimensional than an equally adequate linear
representation. We conclude that even in nonisometric situations the
proposed approach can often be expected to lead to improved representations
of functional data.

%s5.2 #&#
\subsection{Simulation results}\label{sec52}

We simulate functional data from manifolds $\mani_1$--$\mani_3$ as
introduced in the previous subsection, aiming to study two questions. First,
when the functional data lie on a manifold, whether it is isometric or not,
does the proposed functional manifold approach lead to better (more
parsimonious, better interpretable) representations of the data,
compared to
functional principal component analysis? Second, for noisy functional data
that do not exactly lie on a manifold, how much improvement may one
gain by
adding the data-adaptive penalties implemented by P-ISOMAP, as
described in
Section \ref{sec31}?

%f2 #&#
%
\begin{figure}

\includegraphics{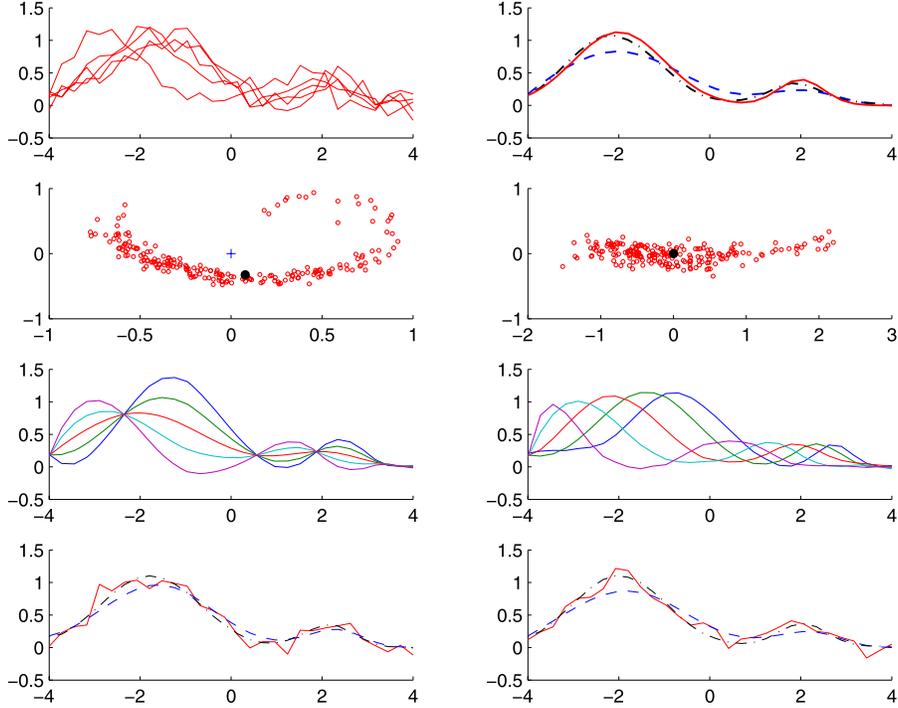}%
\vspace*{-2pt}
\caption{Simulated data for manifold $\mani_1$.
Here and in the following figures, color descriptions
refer to the online version of the paper.
Top left
panel: five randomly selected curves. Top right panel: common shape
function (solid red, corresponds to target mean), estimated manifold
mean $\hat{\mu}^\mani$ (\protect\ref{341}) (dash-dot black) and the
$L^2$ mean
(dashed blue). Second row: scatter plot of second versus first functional
principal component (left) and second versus first functional
manifold component (right), where the bold black dot represents the
manifold mean and the blue cross dot represents the $L^2$ mean. Third
row: estimates of principal component based mode $X_{1,\alpha}$
(\protect\ref{219}) (left) and of manifold mode $X_{1,\alpha}^\mani$
(\protect\ref{2111}) (right) of functional variation for $\alpha=-2, -1,
0, 1, 2$. Bottom row: two randomly selected curves (solid red), with the
corresponding principal component based predictions $\hat{X}^L_i$
(\protect\ref{312}) (dashed blue), and manifold based predictions
$\hat{X}^\mani_i$ (\protect\ref{pred}) (dash-dot black) for
$L=d=2$.}
\label{fig2}
\vspace*{-4pt}
\end{figure}

For these simulations, the actual error-contaminated observations of
the functional trajectories are generated as
$Y_{ij}=X_i(t_{ij})+\varepsilon_{ij}$, $\varepsilon_{ij}\sim
N(0,\sigma^2)$ i.i.d., $i=1,\ldots,n$, $j=1,\ldots,n_i$, where $n=200$,
$t_{ij}$ equally spaced in $[-4,4]$ with $30$ observations per
trajectory, and the noise variance $\sigma^2$ is such that the
signal-to-noise ratio $R$ is $0.1$ or $0.5$. We estimated manifold
means~$\mu^\mani$~(\ref{213}), manifold modes of functional variation
$\Xj$ (\ref{2111}) and obtained\vspace*{1pt} predic-\break tions~%
$\hat{X}_i^\mani$ (\ref {pred}), which were compared with predictions
obtained by functional principal component analysis.

\begin{figure}

\includegraphics{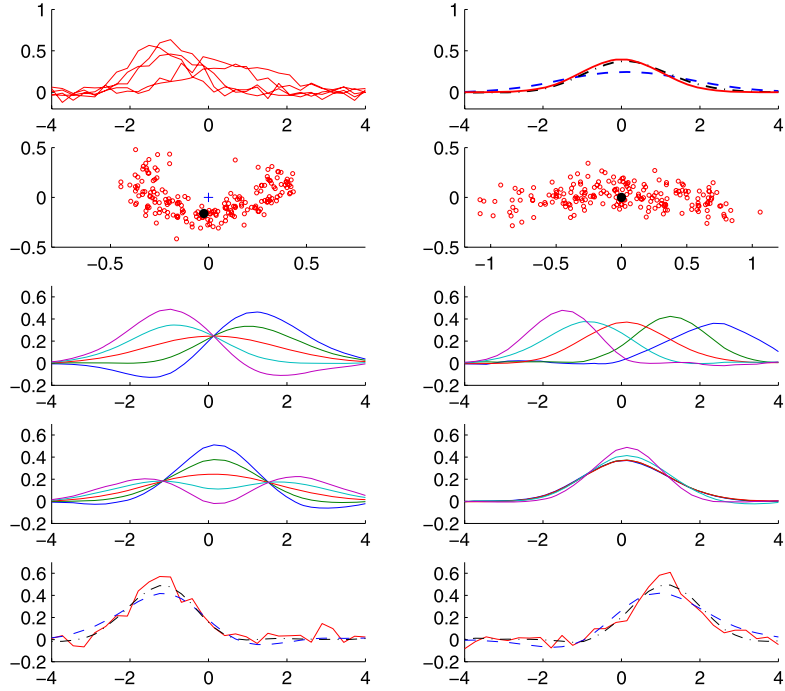}

\caption{Simulated data for manifold $\mani_2$. Top left
panel: five randomly selected curves. Top right panel: standard
Gaussian density (solid red, corresponds to target mean), estimated
manifold mean $\hat{\mu}^\mani$ (\protect\ref{341}) (dash-dot black) and
the $L^2$
mean (dashed blue). Second row: scatter plot of second versus first FPC
(left) and second versus first FMC (right), where the bold black dot
represents the manifold mean and the blue cross represents the~$L^2$ mean.
Third row: estimates of principal component based mode
$X_{1,\alpha}$ (\protect\ref{219}) (left) and of manifold mode
$X_{1,\alpha}^\mani$ (\protect\ref{2111}) (right) of functional variation
for $\alpha=-2, -1, 0, 1, 2$. Fourth row: estimates of
$X_{2,\alpha}$ (left) and of $X_{2,\alpha}^\mani$ (right) for
$\alpha=-2, -1, 0, 1, 2$. Bottom row: two randomly selected curves
(solid red), with the corresponding principal component based predictions
$\hat{X}^L_i$ (\protect\ref{312}) (dash blue), and manifold based predictions
$\hat{X}^\mani_i$ (\protect\ref{pred}) (dash-dot black) for $L=d=3$.}
\label{fig3}
\end{figure}

\begin{figure}

\includegraphics{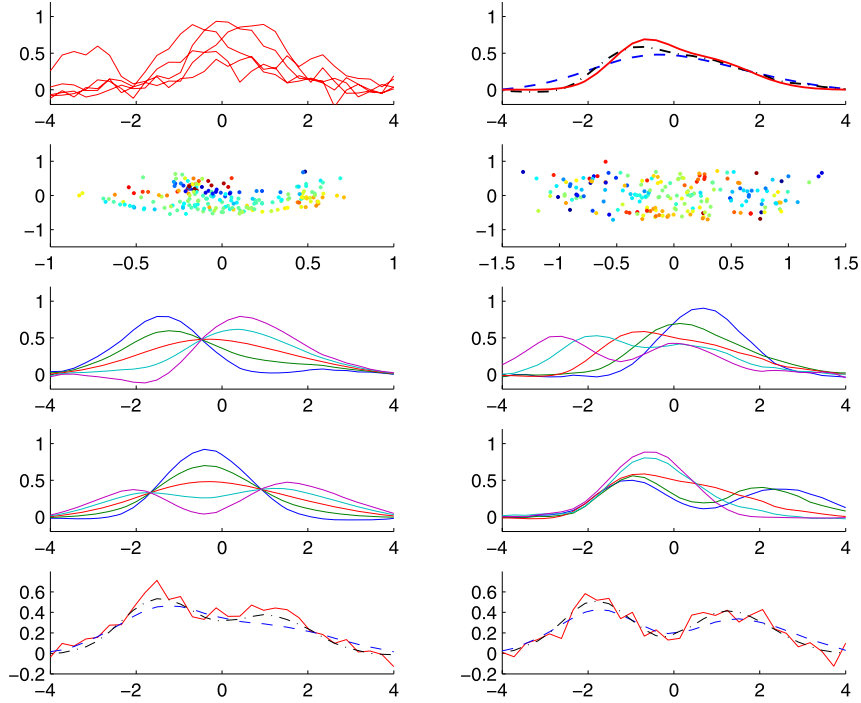}

\caption{Simulated data for manifold $\mani_3$. Top left
panel: five randomly selected curves. Top right panel: curve with no
time shifts (solid red, corresponds to target mean), estimated manifold
mean $\hat{\mu}^\mani$ (\protect\ref{341}) (dash-dot black) and the
$L^2$ mean
(dash blue). Second row: contour scatter plot of second versus first FPC
(left) and second versus first FMC (right), with the colors scaled
from the third FPC or FMC. Third row: estimates of principal
component based mode $X_{1,\alpha}$ (\protect\ref{219}) (left) and of
manifold mode $X_{1,\alpha}^\mani$ (\protect\ref{2111}) (right) of
functional variation for $\alpha=-2, -1, 0, 1, 2$. Fourth row:
estimates of $X_{2,\alpha}$ (left) and of~$X_{2,\alpha}^\mani$ (right)
for $\alpha=-2, -1, 0, 1, 2$. Bottom row: two randomly selected curves
(solid red), with the corresponding
principal component based predictions $\hat{X}^L_i$ (\protect\ref{312})
(dash blue), and manifold based predictions $\hat{X}^\mani_i$
(\protect\ref{pred}) (dash-dot black) for $L=d=3$.}
\label{fig4}
\end{figure}

Results for a simulation run are shown in Figures \ref{fig2}, \ref
{fig3} and \ref{fig4} for manifolds
$\mani_1$--$\mani_3$, respectively. The estimated manifold means are seen
to be
close to the corresponding intrinsic means, that is, the common shape function
for manifold~$\mani_1$, the standard Gaussian density for manifold
$\mani_2$
and the curve with no time shifts ($\alpha=\beta=0$) for manifold
$\mani_3$.
On the other hand, the cross-sectional means are seen to be far away from
these intrinsic means and therefore clearly are not useful as measures of
location for these sets of functions.\vadjust{\goodbreak}

%f3 #&#
%

%f4 #&#
%

The scatter plots of second versus first FPC indicate ``horseshoe'' shapes
for manifolds $\mani_1$ and $\mani_2$. This diagnostic indicates that a
functional manifold approach may be called for. We find that the
location of
the cross-sectional mean (at the origin, due to the zero expectation property
of FPCs) typically lies in a relatively sparse region of the data in these
scatter plots, while the manifold mean falls into a much denser area, which
is another diagnostic feature pointing to an underlying manifold. Complex
two-dimensional surface curvature is observed for manifold $\mani_3$.
Comparing with Figure \ref{fig1}, we find that the manifold modes
represent the
inherent components of functional variation present in the data quite well,
while the established principal component based modes are not
informative in
describing the functional variation. It is also obvious that the proposed
predictions for individual trajectories $X_i$ are more accurate in capturing
amplitudes and locations of peaks.

%
%t2 #&#
%
\begin{table}
\tabcolsep=0pt
\caption{Mean squared prediction errors and relative squared prediction
errors for $\mani_1$--$\mani_3$}\label{table2}
\begin{tabular*}{\tablewidth}{@{\extracolsep{4in minus 4in}}lccd{1.3}ccccd{2.0}d{2.0}d{2.0}d{2.0}c@{}}
\hline
& & &
\multicolumn{5}{c}{\textbf{MSPE with $\bolds{L}$ or $\bolds{d}$}}
& \multicolumn{5}{c@{}}{\textbf{RSPE with $\bolds{L}$ or $\bolds{d}$ (\%)}}\\
[-4pt]
& & &
\multicolumn{5}{c}{\hrulefill} & \multicolumn{5}{c@{}}{\hrulefill}\\
& \multicolumn{1}{c}{$\bolds{R}$} & \textbf{Method}
& \multicolumn{1}{c}{\textbf{1}} & \multicolumn{1}{c}{\textbf{2}}
& \multicolumn{1}{c}{\textbf{3}} & \multicolumn{1}{c}{\textbf{4}}
& \multicolumn{1}{c}{\textbf{5}} & \multicolumn{1}{c}{\textbf{1}}
& \multicolumn{1}{c}{\textbf{2}} & \multicolumn{1}{c}{\textbf{3}}
& \multicolumn{1}{c}{\textbf{4}} & \multicolumn{1}{c@{}}{\textbf{5}} \\
\hline
$\mani_1$ & 0.1&
$\hat{X}_i^L$
& 0.159 & 0.034 & 0.025 & 0.021 & 0.021
& 41 & 10 & 6 & 6 & \hphantom{0}6 \\[2pt]
& & $\hat{X}_i^\mani$ & 0.027 & 0.015 & 0.015 & 0.014 & 0.015
& 7 & 4 & 4 & 4 & \hphantom{0}4 \\
[4pt]
& 0.5 & $\hat{X}_i^L$ & 0.173 & 0.061 & 0.057 & 0.058 & 0.058
& 45 & 16 & 15 & 15 & 15 \\[2pt]
& & $\hat{X}_i^\mani$ & 0.090 & 0.046 & 0.046 & 0.049 & 0.053
& 23 & 12 & 12 & 13 & 14 \\
[4pt]
$\mani_2$ & 0.1 &
$\hat{X}_i^L$
& 0.054 & 0.022 & 0.013 & 0.008 & 0.007
& 44 & 17 & 10 & 7 & \hphantom{0}6 \\[2pt]
& & $\hat{X}_i^\mani$ & 0.022 & 0.009 & 0.007 & 0.006 & 0.006
& 18 & 8 & 5 & 5 & \hphantom{0}5 \\
[4pt] & 0.5& $\hat{X}_i^L$ & 0.055 & 0.025 & 0.019 & 0.018 & 0.018
& 45 & 20 & 16 & 14 & 14 \\[2pt]
& & $\hat{X}_i^\mani$ & 0.030 & 0.017 & 0.015 & 0.014 & 0.014
& 25 & 13 & 12 & 12 & 12 \\
[4pt]
$\mani_3$ & 0.1 & $\hat{X}_i^L$
& 0.148 & 0.059 & 0.031 & 0.023 & 0.020
& 59 & 23 & 13 & 9 & \hphantom{0}7 \\[2pt]
& & $\hat{X}_i^\mani$ & 0.088 & 0.025 & 0.020 & 0.020 & 0.019
& 35 & 10 & 8 & 8 & \hphantom{0}8 \\
[4pt]
& 0.5 & $\hat{X}_i^L$ & 0.154 & 0.071 & 0.053 & 0.048 & 0.048
& 61 & 28 & 21 & 19 & 19 \\[2pt]
& & $\hat{X}_i^\mani$ & 0.124 & 0.059 & 0.047 & 0.045 & 0.044
& 49 & 24 & 19 & 18 & 18 \\
\hline
\end{tabular*}
\vspace*{-3pt}
\end{table}

Leave-one-out predictions of the $X_i$ are calculated using both functional
principal components (\ref{312}), resulting in $\hat{X}_i^L$, as
well as
the proposed new estimates $\hat{X}_i^\mani$ (\ref{pred}). For $\hat
{X}_i^L$,
we estimate the FPCs (\ref{217}) of $X_i$ using all data and then leave
$X_i$ out to obtain $\hat{\mu}$ and $\hat{\phi}_k$; for $\hat
{X}_i^\mani
$, we
estimate $\hat{\psi}(X_i)$ using all data and then leave $X_i$ out in the
local averaging step. Starting with $L=1, d=1$, we increase $L$ and $d$
successively, obtaining the mean squared prediction errors
$\mathrm{MSPE}=\frac{1}{200}\sum_{i=1}^{200} \|X_i-\hat
{X}_i\|^2_{L^2}$, where
$\hat{X}_i=\hat{X}_i^L$ or $\hat{X}_i^\mani$, for $1 \le d=L \le 5$.

The simulation results for manifolds $\mani_1$--$\mani_3$ are shown in
Table \ref{table2}. Generally, the MSPE is reduced by $20\%$ over the
established linear method when using the manifold approach; this
improvement exceeds $50\% $ when~$L$ and~$d$ are small. Another metric
of interest is the relative squared prediction error of the model over
the squared error when using the mean as predictor,
$\mathrm{RSPE}=\frac{\sum_{i=1}^{200}
\|X_i-\hat{X}_i\|^2_{L^2}}{\sum_{i=1}^{200} \|X_i-\bar {X}\|^2_{L^2}}$,
where $\bar{X}=\frac{1}{200}\sum_{i=1}^{200}X_i$, which can be
interpreted as fraction of variance that is left unexplained.
%functional manifold model, the fraction of unexplained variance
%explained
%cannot be attributed to individual components as in the in the linear
%functional principal component model.
In all three simulated manifolds,
RSPE is found to be much larger for the functional principal component
representations, when the same number of components is used. This is because
in the inefficient linear representation higher order functional principal
components carry substantial variation.

To quantify the efficiency of the data-adaptive penalties in the
proposed P-ISOMAP procedure, we also calculated the MSPE using the
unmodified ISOMAP. Parameters for ISOMAP were selected analogously to
the description in Section~\ref{sec33} by cross-validation. Since the
most important comparison is for the case where $d$ equals the
intrinsic dimension, that is, $1$ for $\mani _1$ and $2$ for $\mani_2$
and $\mani_3$, we calculated the ratio of the MSPE of P-ISOMAP over the
MSPE of ISOMAP for these situations (Table \ref {table3}). As
anticipated, P-ISOMAP indeed exhibits increasing benefits for smaller
signal-to-noise ratios.

The influence of the selection of the step size parameter $\varepsilon$
in P-ISOMAP, defined in (\ref{def43}), on mean squared
prediction\vadjust{\goodbreak}
errors is demonstrated in Table~\ref{table4}. Here $d$ is fixed as the
intrinsic dimension ($1$ for $\mani_1$ and 2 for $\mani_2$, $\mani_3$),
while $\delta$ and $h$ are optimized by cross-validation for each
$\varepsilon$. We then select $\varepsilon$ from the median distances
of the $3$rd, $5$th, $8$th, $12$th and $16$th nearest points calculated
over all sample data. From the results in the table, one finds that the
results are not strongly sensitive to the selection of $\varepsilon$,
as long as it is in medium range. A good overall choice is median
distance of $8$th nearest neighbors. When $\varepsilon$ is chosen very
small, some sample points that are not situated close to other sample
points may become separated from the other data, or disconnected
subgroups in the data may emerge, which renders the MSPE for small
$\varepsilon$ inaccurate. In practice, we therefore impose a~lower
bound on $\varepsilon$ to ensure that the fraction of data that are not
connected to other points when connecting through
$\varepsilon$-neighborhoods stays below $5\%$.

%f5 #&#
%
\begin{figure}

\includegraphics{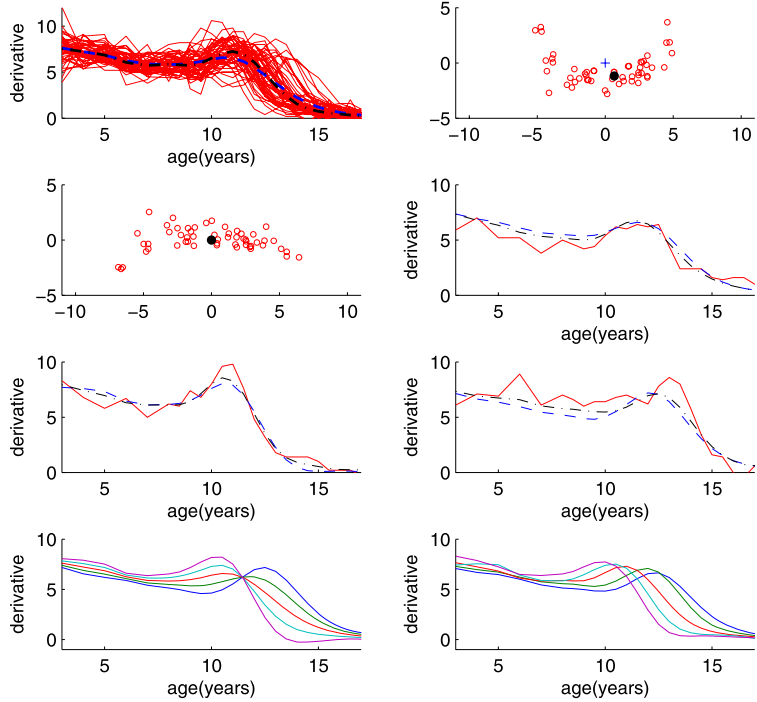}

\caption{Berkeley growth data for girls. Top left\vspace*{1pt} panel:
derivatives with the cross-sectional mean (dash blue) and estimated
manifold mean $\hat{\mu}^\mani$ (\protect\ref{341}) (dash-dot black).
Top right
panel: scatter plot of second versus first FPC, where the bold black
dot represents the manifold mean and the blue cross represents
the cross-sectional mean. Second row left panel: scatter plot of
second versus first FMC, where the bold black dot represents the
manifold mean. Second row right panel and third row panels: three
randomly selected curves (solid red), with the corresponding principal
component based predictions $\hat{X}^L_i$ (\protect\ref{312}) (dash
blue), and
manifold based predictions $\hat{X}^\mani_i$ (\protect\ref{pred})
(dash-dot black)
for $L=d=2$. Bottom panels: estimates of principal component based
mode $X_{1,\alpha}$ (\protect\ref{219}) (left) and of manifold mode
$X_{1,\alpha}^\mani$ (\protect\ref{2111}) (right) of functional variation
for $\alpha=-2, -1, 0, 1, 2$.}
\label{fig5}
\end{figure}

%
%
%%t3 #&#
%%
\begin{table}
\tablewidth=140pt
\caption{Mean squared prediction error ratios
for P-ISOMAP over ISOMAP}\label{table3}
\begin{tabular*}{\tablewidth}{@{\extracolsep{\fill}}lccc@{}}
\hline
$\bolds{R}$ & $\bolds{\mani_1}$ & $\bolds{\mani_2}$ & $\bolds{\mani_3}$ \\
\hline
0.1 & 0.9676 & 0.9679 & 0.9402 \\
0.5 & 0.8121 & 0.8879 & 0.8302 \\
\hline
\end{tabular*}
\end{table}
%
%
%%t4 #&#
%%
\begin{table}[b]
\tablewidth=\textwidth
\caption{Mean squared prediction errors using different
$\varepsilon$ for P-ISOMAP}\label{table4}
\begin{tabular*}{\tablewidth}{@{\extracolsep{\fill}}lccccd{1.3}c@{}}
\hline
& & \multicolumn{5}{c@{}}{$\bolds{\varepsilon}$} \\[-4pt]
& & \multicolumn{5}{c@{}}{\hrulefill} \\
\textbf{Manifold} & $\bolds{R}$ & \textbf{3} & \textbf{5} &
\textbf{8} & \multicolumn{1}{c}{\textbf{12}} & \textbf{16} \\
\hline
$\mani_1$& 0.1 & 0.029 & 0.027 & 0.031 & 0.029 &
0.033\\
& 0.5 & 0.116 & 0.102 & 0.090 & 0.129 & 0.135 \\
[3pt]
$\mani_2$& 0.1 & 0.008 & 0.010 & 0.009 & 0.010 &
0.010 \\
& 0.5 & 0.020 & 0.018 & 0.018 & 0.017 & 0.017 \\
[3pt]
$\mani_3$ & 0.1 & 0.029 & 0.040 & 0.025 & 0.27 &
0.033\\
& 0.5 & 0.052 & 0.059 & 0.065 & 0.059 & 0.066 \\
\hline
\end{tabular*}
\end{table}
%

%s6 #&#
\section{Applications}\label{sec6}

%s6.1 #&#
\subsection{Berkeley growth study}\label{sec61}

In growth studies, one often observes phase variation in the
trajectories. Some subjects reach certain growth stages (such as
puberty in human growth) earlier than others. This leads to
difficulties for the parsimonious modeling of growth
patterns\vspace*{1pt} with linear methods, and more generally for
methods that are based on $L^2$ distance between trajectories.
Accordingly, cross-sectional mean estimation tends to fail in
representing important growth features adequately \cite{kne92,ger05}.
Since phase variation introduces nonlinear features in functional data,
it is of interest to determine whether the analysis of growth data may
benefit from the manifold approach.

%f6 #&#
%
\begin{figure}

\includegraphics{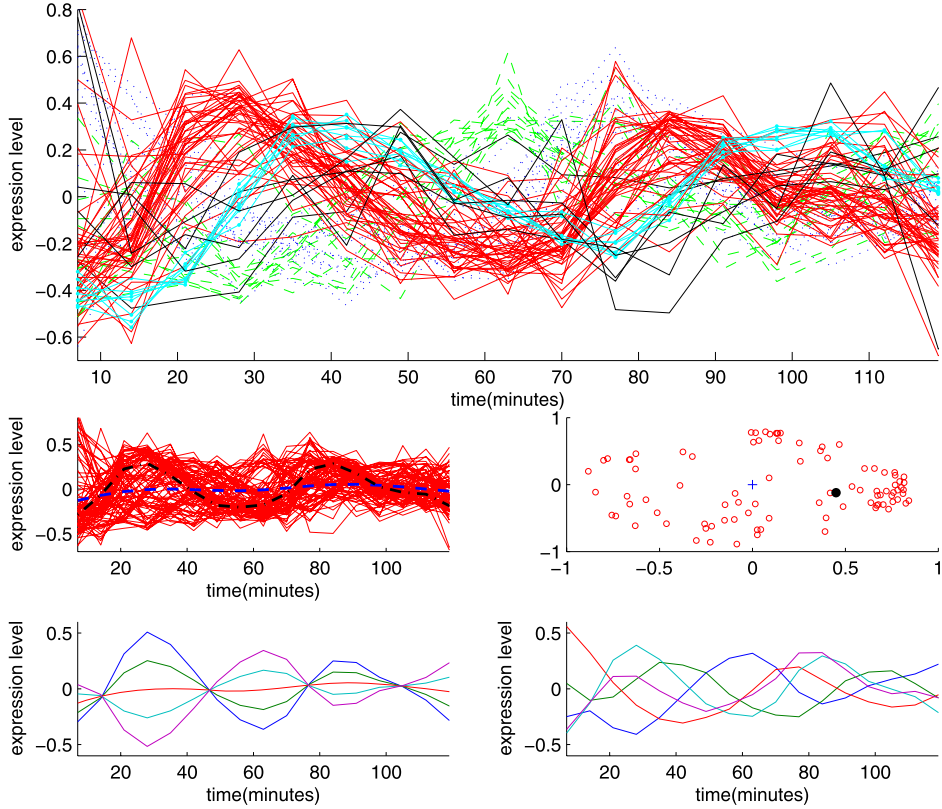}%
\caption{Yeast cell cycle gene expression data. Top panel:
all trajectories in different colors according to cluster
membership: G1 (solid red), S (dash-dot cyan), G2/M (dash green),\vspace*{1pt} M/G1
(dotted blue) and S/G2
(solid black). Middle left panel: estimated manifold
mean~$\hat{\mu}^\mani$ (\protect\ref{341}) (dash-dot black) and
cross-sectional mean
(dash blue). Middle right panel: scatter plot of second versus first FPC,
where the blue cross indicates cross-sectional mean and the bold
black dot indicates manifold mean. Bottom panels: estimates of
principal component based mode $X_{1,\alpha}$ (\protect\ref{219}) (left)
and of manifold mode $X_{1,\alpha}^\mani$ (\protect\ref{2111}) (right) of
functional variation
for $\alpha=-2, -1, 0, 1, 2$.}
\label{fig6}
\vspace*{3pt}
\end{figure}

\begin{figure}

\includegraphics{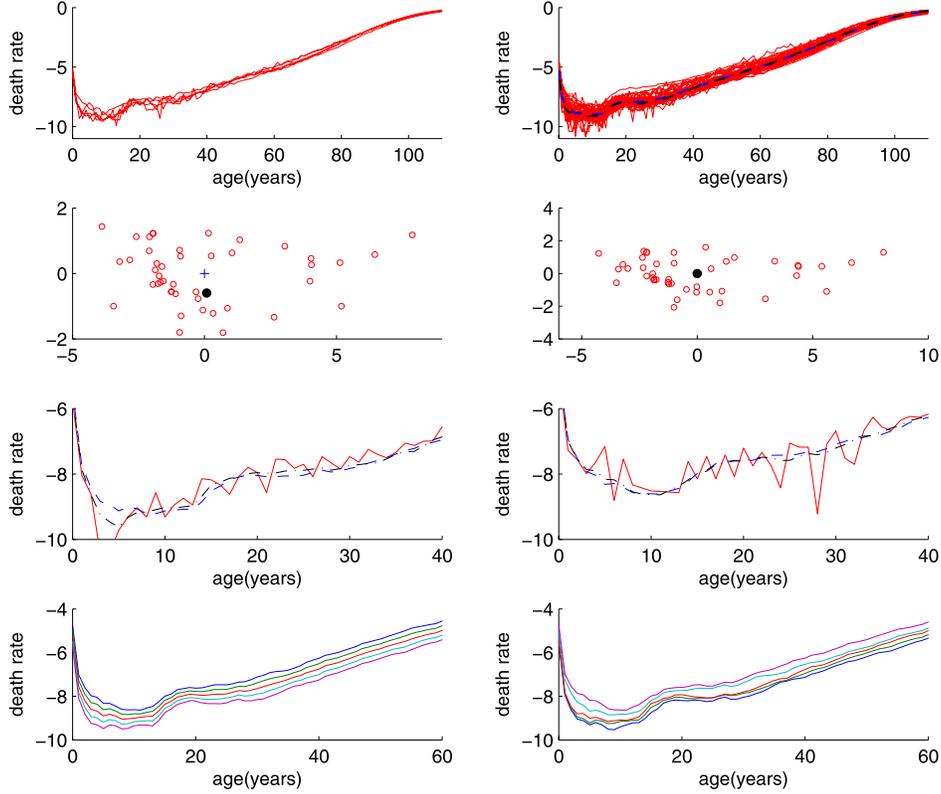}

\caption{Human mortality data. Top left panel: death rates
for five randomly selected countries. Top right panel: estimates of
cross-sectional mean (dash blue) and manifold mean $\hat{\mu}^\mani$
(\protect\ref{314}) (dash{-}dot black). Second row: scatter plots of
second versus
first FPC (left) and second versus first FMC (right), where the blue
cross indicates the cross-sectional mean and the bold black dot
indicates the manifold mean. Third row: two randomly selected curves
(solid red), with the corresponding principal component based predictions
$\hat{X}^L_i$ (\protect\ref{312}) (dash blue), and manifold based predictions
$\hat{X}^\mani_i$ (\protect\ref{pred}) (dash-dot black) for $L=d=3$.
Bottom panels:
estimates of principal component based mode $X_{1,\alpha}$
(\protect\ref{219}) (left) and of manifold mode $X_{1,\alpha}^\mani$
(\protect\ref{2111}) (right) of functional variation for $\alpha=-2, -1,
0, 1, 2$.}
\label{fig7}
\vspace*{3pt}
\end{figure}

We apply the manifold approach to the Berkeley growth data for
females~\cite{tud54}. The data contain height measurements for $54$ girls, with
$31$ measurements taken between the ages of $1$ and $18$ years. Interest
usually focuses on growth velocity \cite{gas84}, which we obtain by
smoothing the first-order difference quotients of the curves. The resulting
growth velocity curves are shown in the top left\vspace*{1pt} panel of Figure \ref{fig5},
together with the cross-sectional mean and the estimated manifold mean
$\hat{\mu}^\mani$ (\ref{341}).
Similarly to Figures \ref{fig2}--\ref{fig4}, the descriptions of
Figures \ref{fig5}--\ref{fig7} refer to the color online versions.
The location of the
cross-sectional mean,
which falls at $(0,0)$, and the location of the estimated manifold mean are
indicated in the scatter plot of second versus first FPC (top right panel),
which displays the ``horseshoe'' pattern described above. This, and the fact
that the cross-sectional mean is away from the main data cloud, point to
inherent nonlinearity in these data.

%f7 #&#
%

Mean squared prediction errors (MSPE) and relative squared prediction errors
(RSPE) for the leave-one-out predictions of $X_i$, as described in
Section~\ref{sec5}, are listed in Table \ref{table5}. The fractions of distance
explained (FDE),
defined in (\ref{331}), for different dimensions $d$ are shown in Table
\ref{table1}. The MSPE of $\hat{X}^L_j$ is minimized at $L=5$, with
$L=2$ already a
quite good choice.\vadjust{\goodbreak}

We find\vspace*{-2pt} that $\hat{X}^\mani_j$ consistently improves upon $\hat
{X}^L_j$, the
fit obtained from functional principal components. Note that we used the
preliminary estimator $\hat{X}^K_i$ in (\ref{pred}) with $K=4$, applying
criterion (\ref{313}). The FDE criterion indicates that these data
can be
well described by a one-dimensional manifold. The middle three panels of
Figure \ref{fig5} include three randomly selected curves, along with
the predictions
$\hat{X}^L_i$ and $\hat{X}^\mani_i$ using $L=d=2$. The two bottom
panels of
Figure \ref{fig5} illustrate the comparison of estimated manifold
modes of
functional variation with the principal component based modes. The manifold
modes are clearly more useful and adequately reflect the time-warping feature
of these data. The first manifold mode specifically suggests that for girls,
a puberty growth peak at a late age, especially after age $12$, tends
to have
a smaller amplitude; this is in line with auxological knowledge.
Overall, the
manifold mode is seen to provide a clearer and much more adequate description
of the longitudinal dynamics of these data.

%s6.2 #&#
\subsection{Yeast cell cycle gene expression}\label{sec62}

Temporal expression curves for yeast cell cycle related genes were obtained
by \cite{spe98}. There are 6,178 genes in total, where each gene expression
time-course consists of $18$ data points, measured every $7$ minutes between
$0$ and $119$ minutes. Groups of genes are thought to be coexpressed
coherently across different time periods, according to the role played
by the
genes in the time progression of the cell cycle. The dynamics of the gene
expression levels are complex. Temporal regularization of gene
expression is
a characteristic of gene function, suggesting models that incorporate
time-warping \cite{len06,tan09}.

The data we study consist of $90$ genes that have been identified by
biological methods \cite{spe98}. Of these genes, $44$ are thought to be
related to G1 phase regulation of the yeast cell cycle and $46$ to non-G1
phase regulation (S, S/G2, G2/M and M/G1 phases). Time courses of gene
expression (top panel of Figure \ref{fig6}) for these clusters reveal
two peaks for
the G1 (solid red) and S (dash-dot cyan) groups, and one peak for G2/M
(dash green) and M/G1
(dotted blue) groups, while the trajectories for the S/G2 (solid black)
group are highly
variable with no obvious peak.

%t5 #&#
%
\begin{table}
\tabcolsep=0pt
\caption{Mean squared prediction errors and relative squared prediction
errors for growth, yeast and mortality data}\label{table5}
\begin{tabular*}{\tablewidth}{@{\extracolsep{4in minus 4in}}lcd{2.3}d{2.3}d{2.3}d{2.3}d{2.3}ccccc@{}}
\hline
& & \multicolumn{5}{c}{\textbf{MSPE with $\bolds{L}$ or $\bolds{d}$}} &
\multicolumn{5}{c@{}}{\textbf{RSPE with $\bolds{L}$ or $\bolds{d}$ (\%)}} \\[-4pt]
& & \multicolumn{5}{c}{\hrulefill} & \multicolumn{5}{c@{}}{\hrulefill}\\
\textbf{Data} & \textbf{Method} &
\multicolumn{1}{c}{\textbf{1}} & \multicolumn{1}{c}{\textbf{2}}
& \multicolumn{1}{c}{\textbf{3}} & \multicolumn{1}{c}{\textbf{4}}
& \multicolumn{1}{c}{\textbf{5}} & \multicolumn{1}{c}{\textbf{1}}
& \multicolumn{1}{c}{\textbf{2}} & \multicolumn{1}{c}{\textbf{3}}
& \multicolumn{1}{c}{\textbf{4}} & \multicolumn{1}{c@{}}{\textbf{5}} \\
\hline
Growth & $\hat{X}_i^L$ & 17.1 & 12.9 & 13.8 &
13.7 & 12.6 & 62 & 47 & 50 & 50 & 46 \\[2pt]
& $\hat{X}_i^\mani$ & 10.7 & 9.46 & 9.06 & 9.21 & 9.08 & 39
& 34 & 33 & 33 & 33 \\
[4pt]
Yeast & $\hat{X}_i^L$ & 0.639 & 0.382 & 0.257 &
0.205 & 0.203 & 67 & 40 & 27 & 22 & 21 \\[2pt]
& $\hat{X}_i^\mani$ & 0.468 & 0.278 & 0.231 & 0.210 & 0.206 &
49 & 29 & 24 & 22 & 22 \\
[4pt]
Mortality & $\hat{X}_i^L$ & 7.38 & 6.34 & 5.44 &
5.48 & 5.21 & 54 & 47 & 40 & 40 & 38 \\[2pt]
& $\hat{X}_i^\mani$ & 6.77 & 5.64 & 5.40 & 5.26 & 4.98 &
50 & 41 & 40 & 39 & 37 \\
\hline
\end{tabular*}
\end{table}

The proposed manifold analysis was applied to this set of $90$ genes. The
estimated manifold mean $\hat{\mu}^\mani$ (\ref{341}) (middle left
panel of
Figure \ref{fig6}) is seen to fall within the G1 group (solid red in
the top panel). In
contrast, the cross-sectional mean is almost flat and does not reflect useful
information about these data. We also calculated the MSPE (Table \ref
{table5}) of
$\hat{X}^L_i$ (\ref{312}) and~$\hat{X}^\mani_i$ (\ref{pred}), using
preliminary estimators $\hat{X}^K_i$ with $K=4$ in (\ref{pred}). The
manifold-based predictions are seen to be much better for $d=1$ and $2$,
while they become more similar in performance to $\hat{X}^L_i$ when $d$
increases.

In the two bottom panels of Figure \ref{fig6}, we display the
estimated manifold
mode (right) and the principal component based mode of functional variation
(left). The latter is found to be deceptive, as it indicates amplitude
variation around a few fixed ``knots,'' while the first manifold mode clearly
illustrates the actual temporal variation in the data, which is mainly caused
by phase shifts. Each of the five groups, except the S/G2 group (solid
black), is
well represented by the variation across this manifold mode.

%s6.3 #&#
\subsection{Human mortality across countries}\label{sec63}

The death rates derived from current lifetable cohorts for $44$
countries in
the year $2000$, recorded for each age ranging from $0$ to $110$, have been
collected and are as described in \url{http://www.lifetable.de/}. Death
rates are
widely used for descriptive and analytical purposes in public health, and
cross-country comparisons are of particular interest here.

We view log-transformed annual death rates as noisy measurements of
underlying smooth trajectories. Five sample trajectories are shown in
the top
left panel of Figure \ref{fig7}. The mortality trajectories are densely
sampled, but
the annual rates are quite noisy. We presmoothed\vspace*{1pt} this data,
following~(\ref{311}). The resulting MSPEs for $\hat{X}^\mani_i$ (\ref
{pred}) and
$\hat{X}^L_i$ (\ref{312}) are in Table~\ref{table5}. Manifold-based
prediction is
seen to perform better than linear principal component based prediction,
regardless of the choice of dimension. This is also illustrated by the
panels in the third row of Figure \ref{fig7}, where
predicted trajectories are obtained for $L=d=3$. For these data, the
estimated manifold mean~$\hat{\mu}^\mani$ (\ref{341}) does not differ
dramatically from the cross-sectional mean (top right panel and second row
left panel). However, the first manifold mode of variation (bottom right
panel) indicates that countries with overall lower death rates, or more
specifically, with death rates below the mean curve (solid red),
exhibit less
variation than those with death rates above the mean, especially for ages
from $0$ to $40$. This finding is in line with the skewness that is apparent
in the scatter plots, but is not seen in the principal component based mode
(bottom left panel). The observed gains in prediction error for the manifold
approach provide evidence that substantial nonlinearity is present in these
data.

%s7 #&#
\section{Discussion}\label{sec7}

While the proposed functional manifold implementations were running
relatively fast on a linux server, observing that the computational
complexity of classical MDS is of the order $O(n^3)$, computational
difficulties may arise for truly large sample sizes $n$. In such situations,
one might consider to base the proposed methods on landmark MDS
\cite{sil03}, where one employs landmarks to significantly reduce the
computational complexity.\looseness=-1\vadjust{\goodbreak}

The proposed method relies on two major assumptions: the isometry of the
underlying functional manifold and that the target manifold is close or
identical to the manifold identified by ISOMAP at the sample points. As for
the isometry assumption, even if it is violated, the proposed method proves
to be beneficial, as it often will provide for a much sparser representation
of functional data in comparison with linear methods in cases where the
underlying manifold is nonlinear, even if this manifold is not isometric.
This is discussed in detail in Section \ref{sec51} and borne out by
simulations. As
for the closeness of the ISOMAP solution to the true manifold at the sample
points, this assumption and its underlying justification pertains to ISOMAP
for vector data as proposed in \cite{ten00}.

Starting from the simplifying assumption that the ISOMAP identified manifold
and the target manifold are essentially identical at the sample points, we
proceed to extend the estimation of the manifold function to the entire space
of interest. We note that such simplifying assumptions are often beneficial
when deploying complex statistical methodology, as even when the assumptions
are not completely satisfied, the resulting methodology may turn out to be
more efficient than existing methods.

Overall, we find that the proposed manifold mean and manifold modes of
functional variation provide useful representations that are
competitive with
and often superior over classical linear representations for functional
data. The proposed functional manifold representations thus complement the
established linear representations, notably the Karhunen--Lo\`{e}ve
representation, and in many instances provide more efficient models with
better interpretations.

\begin{appendix}\label{app}

%s8 #&#
\section*{Appendix: Assumptions}

\begin{longlist}[(A1)]
\item[(A1)] The bandwidths $h_{\mu}$, $h_v$, $h_G$ for estimating
$\mu(t)$, $\sigma^2$, $G(t,s)$ in Section~\ref{sec31} satisfy: $h_\mu
\rightarrow0$, $nh_\mu^4 \rightarrow\infty$ and $nh_\mu^6 <
\infty$; $h_G \rightarrow0$, $nh_G^6 \rightarrow\infty$ and $nh_G^8
< \infty$; $h_V \rightarrow0$, $nh_V^4 \rightarrow\infty$ and
$nh_V^6 < \infty$.
\item[(A2)] The smoothing kernels $\kappa_{\mu}$ for the mean function
$\mu$ and $\kappa_G$ for the covariance function $G$ in Section
\ref{sec31}
are absolutely integrable, that is, $\int|\kappa_{\mu}(t)|\,dt <\infty$
and $\int\!\!\int|\kappa_G(t,s)|\,dt\,ds <\infty$.
\item[(A3)] For $\dij=t_{ij}-t_{i,j-1}$ and $\tau^*=\max_{i,j}\dij
$, it holds
that $\tau^*=O_p(r^2_n)$, where
$r_n=\max\{\frac{1}{\sqrt{n}h^2_G},\frac{1}{\sqrt{n}h_{\mu}},
\frac{1}{\sqrt{n}h_V}\}$.
\item[(A4)] The eigenvalues of the covariance function $G(t,s)$ satisfy
$\lambda_k<C_1 k^{-\alpha_1}$ for some constants $C_1<\infty$,
$\alpha_1>1$, and if $\lambda_k>0$, then $\lambda_k-\lambda_{k+1}>C_2
k^{-\alpha_2}$ for some constants $C_2>0$ and $\alpha_2>0$.
\item[(A5)] For any $X \in\mani$, $X$ is differentiable and
$\|X\|_{\infty}=O_p(1)$, $\|X'\|_{\infty}=O_p(1)$. The covariance
function $G(t,s)$ is twice differentiable in both $t$ and~$s$, and
$\sup_{t,s\in\mathcal{T}}|G(t,s)|<C_3$, $\sup_{t,s\in
\mathcal{T}}|\frac{\partial^2G(t,s)}{\partial t\,\partial s}| < C_4$
for some constants $C_3,\break C_4<\infty$.
\item[(B1)] The estimates $\hat{\psi}$ of $\psi$ converge uniformly on
the sample space, that is, $\E v_n\rightarrow0$ for
$v_n={\sup_{i=1,\ldots,n}}\|\hat{\psi}(X_i)-\psi(X_i)\|$.
\item[(B2)] Each component of the $d$-vector $\psi(X)$ has a finite
fourth moment, and its covariance matrix is positive definite.
\item[(C1)] The $d$-vector $\psi(X)$ admits a density function $f$, which
is twice differentiable with continuous partial derivatives and
uniformly bounded Hessian matrix.
\item[(C2)] The $d$-dimensional nonnegative kernel $\kappa$ satisfies
$\int\kappa(\bu)\,d\bu=1$,\break $\kappa(\bu)=\kappa(-\bu)$,
$\operatorname{det}
(\int\kappa(\bu)\bu\bu^T \,d\bu)<\infty$, $\int\kappa^2(\bu
)\,d\bu<\infty$,
and is Lipschitz continuous with compact support, $\{\bu\in\R
^d\dvtx\|\bu
\|\leq1\}$.
\item[(C3)] The map $\psi^{-1}\dvtx\R^d\rightarrow L^2$ is twice Fr\'{e}chet
differentiable, that is, there exist bounded linear operators
$A^1_{\mathbf{u}}\dvtx\R^d\rightarrow L^2$, $A^2_{{\mf u}}\dvtx\R
^d\times\R^d\rightarrow
L^2$ such that
\begin{eqnarray*}
\lim_{{\mf u}_1\rightarrow{\mf
0}}\frac{\|\psi^{-1}({\mf u}+{\mf u}_1)-\psi^{-1}({\mf u})-A^1_{{\mf
u}} ({\mf u}_1)\|_{L^2}}{\|{\mf u}_1\|}&=&0,
\\
\lim_{{\mf u}_2\rightarrow{\mf0}}\frac{\|A^1_{{\mf u}+{\mf u}_2}({\mf
u}_1)-A^1_{{\mf u}+{\mf u}_2}({\mf u}_1)- A^2_{{\mf u}}({\mf u}_1,{\mf
u}_2)\|_{L^2}}{\|{\mf u}_2\|}&=&0
\end{eqnarray*}
for all ${\mf u}, {\mf u}_1, {\mf u}_2 \in\R^d$. In addition, $\frac{\|
A^2_{{\mf u}}({\mf u}_1,{\mf u}_2)\|_{L^2}}{\|{\mf u}_1\|\cdot\|{\mf
u}_2\|}$ is continuous
and uniformly bounded w.r.t. ${\mf u}$.
\end{longlist}
\end{appendix}

\section*{Acknowledgments}

We are grateful for helpful comments on previous versions of this paper
from two referees and an Associate Editor.

\begin{supplement}[id=suppA]
\stitle{Supplement to ``Nonlinear manifold
representations for functional data''}
\slink[doi]{10.1214/11-AOS936SUPP} %[doi,text={...}] - jei reikia
%suskaldyti doi
\sdatatype{.pdf}
\sfilename{aos936\_supp.pdf}
\sdescription{An online supplementary file contains the detailed
proofs for Propositions \ref{prop1}--\ref{prop4}, Theorem \ref{thm1}
and Corollary \ref{Cor1}.
These proofs make use of material in references
\cite{hal07,han08,kat66,mar79}.}
\end{supplement}

% imsref loaded by lrinkeviciute, 2011-12-20 15:21:21
% imsref loaded by lrinkeviciute, 2011-12-20 15:24:27
%

\printaddresses


\begin{thebibliography}{35}
% BibTex style file: ims.bst, 2011-05-30
% Default style options (sort=0,type=number).
% Used options (sort=1,type=number).

%b1 #&#
\bibitem{ash75}
%
\begin{bbook}[mr]
\bauthor{\bsnm{Ash},~\bfnm{Robert~B.}\binits{R.~B.}} \AND
\bauthor{\bsnm{Gardner},~\bfnm{Melvin~F.}\binits{M.~F.}}
(\byear{1975}).
\btitle{Topics in Stochastic Processes}.
\bseries{Probability and Mathematical Statistics}
\bvolume{27}.
\bpublisher{Academic Press},
\baddress{New York}.
\bid{mr={0448463}}
\bptok{imsref}%
\end{bbook}
%
\endbibitem

%b2 #&#
\bibitem{bel03}
%
\begin{barticle}[auto:STB|2011/12/15|13:36:40]
\bauthor{\bsnm{Belkin},~\bfnm{M.}\binits{M.}} \AND
\bauthor{\bsnm{Niyogi},~\bfnm{P.}\binits{P.}}
(\byear{2003}).
\btitle{Laplacian eigenmaps for dimensionality reduction and data
representation}.
\bjournal{Neural Comput.}
\bvolume{15}
\bpages{1373--1396}.
\bptok{imsref}%
\end{barticle}
%
\endbibitem

%b3 #&#
\bibitem{bic07}
%
\begin{bincollection}[mr]
\bauthor{\bsnm{Bickel},~\bfnm{Peter~J.}\binits{P.~J.}} \AND
\bauthor{\bsnm{Li},~\bfnm{Bo}\binits{B.}}
(\byear{2007}).
\btitle{Local polynomial regression on unknown manifolds}.
In \bbooktitle{Complex Datasets and Inverse Problems}.
\bseries{IMS Lecture Notes Monogr. Ser.}
\bvolume{54}
\bpages{177--186}.
\bpublisher{IMS}, \baddress{Beachwood, OH}.
\bid{doi={10.1214/074921707000000148}, mr={2459188}}
\bptok{imsref}%
\end{bincollection}
%
\endbibitem

%b4 #&#
\bibitem{cas86}
%
\begin{barticle}[auto:STB|2011/12/15|13:36:40]
\bauthor{\bsnm{Castro},~\bfnm{P.~E.}\binits{P.~E.}},
\bauthor{\bsnm{Lawton},~\bfnm{W.~H.}\binits{W.~H.}} \AND
\bauthor{\bsnm{Sylvestre},~\bfnm{E.~A.}\binits{E.~A.}}
(\byear{1986}).
\btitle{Principal modes of variation for processes with continuous sample
curves}.
\bjournal{Technometrics}
\bvolume{28}
\bpages{329--337}.
\bptok{imsref}%
\end{barticle}
%
\endbibitem

%b5 #&#
\bibitem{supp}
%
\begin{bmisc}[auto:STB|2011/12/15|13:36:40]
\bauthor{\bsnm{Chen},~\bfnm{D.}\binits{D.}} \AND
\bauthor{\bsnm{M{\"u}ller},~\bfnm{H.~G.}\binits{H.~G.}}
(\byear{2011}).
\bhowpublished{Supplement to ``Nonlinear manifold representations for
functional data.''
\href{http://dx.doi.org/10.1214/11-AOS936SUPP}{DOI:10.1214/11-AOS936SUPP}.}
\bptok{imsref}%
\end{bmisc}
%
\endbibitem

%b6 #&#
\bibitem{cox01}
%
\begin{bbook}[mr]
\bauthor{\bsnm{Cox},~\bfnm{Trevor~F.}\binits{T.~F.}} \AND
\bauthor{\bsnm{Cox},~\bfnm{Michael A.~A.}\binits{M.~A.~A.}}
(\byear{2001}).
\btitle{Multidimensional Scaling}.
\bpublisher{Chapman and Hall}, \baddress{London}.
\bptnote{check year}%
\bptok{imsref}%
\end{bbook}
%
\endbibitem

%b7 #&#
\bibitem{sil03}
%
\begin{barticle}[auto:STB|2011/12/15|13:36:40]
\bauthor{\bsnm{De~Silva},~\bfnm{V.}\binits{V.}} \AND
\bauthor{\bsnm{Tenenbaum},~\bfnm{J.~B.}\binits{J.~B.}}
(\byear{2003}).
\btitle{Global versus local methods in nonlinear dimensionality reduction}.
\bjournal{Adv. Neural Inf. Process. Syst.}
\bvolume{15}
\bpages{721--728}.
\bptok{imsref}%
\end{barticle}
%
\endbibitem

%b8 #&#
\bibitem{doc92}
%
\begin{bbook}[mr]
\bauthor{\bparticle{do}
\bsnm{Carmo},~\bfnm{Manfredo~Perdig{\~a}o}\binits{M.~P.}}
(\byear{1992}).
\btitle{Riemannian Geometry}.
\bpublisher{Birkh\"auser}, \baddress{Boston, MA}.
%Flaherty}.
\bid{mr={1138207}}
\bptok{imsref}%
\end{bbook}
%
\endbibitem

%b9 #&#
\bibitem{don03}
%
\begin{barticle}[mr]
\bauthor{\bsnm{Donoho},~\bfnm{David~L.}\binits{D.~L.}} \AND
\bauthor{\bsnm{Grimes},~\bfnm{Carrie}\binits{C.}}
(\byear{2003}).
\btitle{Hessian eigenmaps: Locally linear embedding techniques for
high-dimensional data}.
\bjournal{Proc. Natl. Acad. Sci. USA}
\bvolume{100}
\bpages{5591--5596 (electronic)}.
\bid{doi={10.1073/pnas.1031596100}, issn={1091-6490}, mr={1981019}}
\bptok{imsref}%
\end{barticle}
%
\endbibitem

%b10 #&#
\bibitem{don05}
%
\begin{barticle}[mr]
\bauthor{\bsnm{Donoho},~\bfnm{David~L.}\binits{D.~L.}} \AND
\bauthor{\bsnm{Grimes},~\bfnm{Carrie}\binits{C.}}
(\byear{2005}).
\btitle{Image manifolds which are isometric to {E}uclidean space}.
\bjournal{J. Math. Imaging Vision}
\bvolume{23}
\bpages{5--24}.
\bid{doi={10.1007/s10851-005-4965-4}, issn={0924-9907}, mr={2208906}}
\bptok{imsref}%
\end{barticle}
%
\endbibitem

%b11 #&#
\bibitem{gas84}
%
\begin{barticle}[mr]
\bauthor{\bsnm{Gasser},~\bfnm{Theo}\binits{T.}},
\bauthor{\bsnm{M{\"u}ller},~\bfnm{Hans-Georg}\binits{H.-G.}},
\bauthor{\bsnm{K{\"o}hler},~\bfnm{Walter}\binits{W.}},
\bauthor{\bsnm{Molinari},~\bfnm{Luciano}\binits{L.}} \AND
\bauthor{\bsnm{Prader},~\bfnm{Andrea}\binits{A.}}
(\byear{1984}).
\btitle{Nonparametric regression analysis of growth curves}.
\bjournal{Ann. Statist.}
\bvolume{12}
\bpages{210--229}.
\bid{doi={10.1214/aos/1176346402}, issn={0090-5364}, mr={0733509}}
\bptok{imsref}%
\end{barticle}
%
\endbibitem

%b12 #&#
\bibitem{ger04}
%
\begin{barticle}[mr]
\bauthor{\bsnm{Gervini},~\bfnm{Daniel}\binits{D.}} \AND
\bauthor{\bsnm{Gasser},~\bfnm{Theo}\binits{T.}}
(\byear{2004}).
\btitle{Self-modelling warping functions}.
\bjournal{J. R. Stat. Soc. Ser. B Stat. Methodol.}
\bvolume{66}
\bpages{959--971}.
\bid{doi={10.1111/j.1467-9868.2004.B5582.x}, issn={1369-7412}, mr={2102475}}
\bptok{imsref}%
\end{barticle}
%
\endbibitem

%b13 #&#
\bibitem{ger05}
%
\begin{barticle}[mr]
\bauthor{\bsnm{Gervini},~\bfnm{Daniel}\binits{D.}} \AND
\bauthor{\bsnm{Gasser},~\bfnm{Theo}\binits{T.}}
(\byear{2005}).
\btitle{Nonparametric maximum likelihood estimation of the structural
mean of a
sample of curves}.
\bjournal{Biometrika}
\bvolume{92}
\bpages{801--820}.
\bid{doi={10.1093/biomet/92.4.801}, issn={0006-3444}, mr={2234187}}
\bptok{imsref}%
\end{barticle}
%
\endbibitem

%b14 #&#
\bibitem{gre50}
%
\begin{barticle}[mr]
\bauthor{\bsnm{Grenander},~\bfnm{Ulf}\binits{U.}}
(\byear{1950}).
\btitle{Stochastic processes and statistical inference}.
\bjournal{Ark. Mat.}
\bvolume{1}
\bpages{195--277}.
\bid{issn={0004-2080}, mr={0039202}}
\bptok{imsref}%
\end{barticle}
%
\endbibitem

%b15 #&#
\bibitem{hal07}
%
\begin{barticle}[mr]
\bauthor{\bsnm{Hall},~\bfnm{Peter}\binits{P.}} \AND
\bauthor{\bsnm{Horowitz},~\bfnm{Joel~L.}\binits{J.~L.}}
(\byear{2007}).
\btitle{Methodology and convergence rates for functional linear regression}.
\bjournal{Ann. Statist.}
\bvolume{35}
\bpages{70--91}.
\bid{doi={10.1214/009053606000000957}, issn={0090-5364}, mr={2332269}}
\bptok{imsref}%
\end{barticle}
%
\endbibitem

%b16 #&#
\bibitem{hel01}
%
\begin{bbook}[mr]
\bauthor{\bsnm{Helgason},~\bfnm{Sigurdur}\binits{S.}}
(\byear{2001}).
\btitle{Differential Geometry, {L}ie Groups, and Symmetric Spaces}.
\bseries{Graduate Studies in Mathematics}
\bvolume{34}.
\bpublisher{Amer. Math. Soc.}, \baddress{Providence, RI}.
\bid{mr={1834454}}
\bptok{imsref}%
\end{bbook}
%
\endbibitem

%b17 #&#
\bibitem{huc11}
%
\begin{barticle}[auto:STB|2011/12/15|13:36:40]
\bauthor{\bsnm{Huckemann},~\bfnm{S.}\binits{S.}}
(\byear{2011}).
\btitle{Inference on 3d Procrustes means: Tree bole growth, rank deficient
diffusion tensors and perturbation models}.
\bjournal{Scand. J. Stat.}
\bvolume{38}
\bpages{1467--9469}.
\bptok{imsref}%
\end{barticle}
%
\endbibitem

%b18 #&#
\bibitem{ize07}
%
\begin{barticle}[mr]
\bauthor{\bsnm{Izem},~\bfnm{Rima}\binits{R.}} \AND
\bauthor{\bsnm{Marron},~\bfnm{J.~S.}\binits{J.~S.}}
(\byear{2007}).
\btitle{Analysis of nonlinear modes of variation for functional data}.
\bjournal{Electron. J. Stat.}
\bvolume{1}
\bpages{641--676}.
\bid{doi={10.1214/07-EJS080}, issn={1935-7524}, mr={2369029}}
\bptok{imsref}%
\end{barticle}
%
\endbibitem

%b19 #&#
\bibitem{jon92}
%
\begin{barticle}[auto:STB|2011/12/15|13:36:40]
\bauthor{\bsnm{Jones},~\bfnm{M.~C.}\binits{M.~C.}} \AND
\bauthor{\bsnm{Rice},~\bfnm{J.~A.}\binits{J.~A.}}
(\byear{1992}).
\btitle{Displaying the important features of large collections of similar
curves}.
\bjournal{Amer. Statist.}
\bvolume{46}
\bpages{140--145}.
\bptok{imsref}%
\end{barticle}
%
\endbibitem

%b20 #&#
\bibitem{kat66}
%
\begin{bbook}[auto:STB|2011/12/15|13:36:40]
\bauthor{\bsnm{Kato},~\bfnm{T.}\binits{T.}}
(\byear{1966}).
\btitle{Perturbation Theory for Linear Operators}.
\bpublisher{Springer}, \baddress{New York}.
\bid{mr={0203473}}
\bptok{imsref}%
\end{bbook}
%
\endbibitem

%b21 #&#
\bibitem{ken99}
%
\begin{bbook}[mr]
\bauthor{\bsnm{Kendall},~\bfnm{D.~G.}\binits{D.~G.}},
\bauthor{\bsnm{Barden},~\bfnm{D.}\binits{D.}},
\bauthor{\bsnm{Carne},~\bfnm{T.~K.}\binits{T.~K.}} \AND
\bauthor{\bsnm{Le},~\bfnm{H.}\binits{H.}}
(\byear{1999}).
\btitle{Shape and Shape Theory}.
\bpublisher{Wiley}, \baddress{Chichester}.
\bid{doi={10.1002/9780470317006}, mr={1891212}}
\bptok{imsref}%
\end{bbook}
%
\endbibitem

%b22 #&#
\bibitem{kne92}
%
\begin{barticle}[mr]
\bauthor{\bsnm{Kneip},~\bfnm{Alois}\binits{A.}} \AND
\bauthor{\bsnm{Gasser},~\bfnm{Theo}\binits{T.}}
(\byear{1992}).
\btitle{Statistical tools to analyze data representing a sample of curves}.
\bjournal{Ann. Statist.}
\bvolume{20}
\bpages{1266--1305}.
\bid{doi={10.1214/aos/1176348769}, issn={0090-5364}, mr={1186250}}
\bptok{imsref}%
\end{barticle}
%
\endbibitem

%b23 #&#
\bibitem{kne01}
%
\begin{barticle}[mr]
\bauthor{\bsnm{Kneip},~\bfnm{Alois}\binits{A.}} \AND
\bauthor{\bsnm{Utikal},~\bfnm{Klaus~J.}\binits{K.~J.}}
(\byear{2001}).
\btitle{Inference for density families using functional principal component
analysis}.
\bjournal{J. Amer. Statist. Assoc.}
\bvolume{96}
\bpages{519--531}.
\bptnote{check related}%
\bptok{imsref}%
\end{barticle}
%
\endbibitem

%b24 #&#
\bibitem{len06}
%
\begin{barticle}[auto:STB|2011/12/15|13:36:40]
\bauthor{\bsnm{Leng},~\bfnm{X.}\binits{X.}} \AND
\bauthor{\bsnm{M{\"u}ller},~\bfnm{H.~G.}\binits{H.~G.}}
(\byear{2006}).
\btitle{Time ordering of gene co-expression}.
\bjournal{Biostatistics}
\bvolume{7}
\bpages{569--584}.
\bptok{imsref}%
\end{barticle}
%
\endbibitem

%b25 #&#
\bibitem{mar79}
%
\begin{bbook}[mr]
\bauthor{\bsnm{Mardia},~\bfnm{Kantilal~Varichand}\binits{K.~V.}},
\bauthor{\bsnm{Kent},~\bfnm{John~T.}\binits{J.~T.}} \AND
\bauthor{\bsnm{Bibby},~\bfnm{John~M.}\binits{J.~M.}}
(\byear{1979}).
\btitle{Multivariate Analysis}.
\bpublisher{Academic Press},
\baddress{London}.
\bid{mr={0560319}}
\bptok{imsref}%
\end{bbook}
%
\endbibitem

%b26 #&#
\bibitem{han08}
%
\begin{barticle}[mr]
\bauthor{\bsnm{M{\"u}ller},~\bfnm{Hans-Georg}\binits{H.-G.}} \AND
\bauthor{\bsnm{Yao},~\bfnm{Fang}\binits{F.}}
(\byear{2008}).
\btitle{Functional additive models}.
\bjournal{J. Amer. Statist. Assoc.}
\bvolume{103}
\bpages{1534--1544}.
\bid{doi={10.1198/016214508000000751}, issn={0162-1459}, mr={2504202}}
\bptok{imsref}%
\end{barticle}
%
\endbibitem

%b27 #&#
\bibitem{rie90}
%
\begin{bbook}[mr]
\bauthor{\bsnm{Riesz},~\bfnm{Frigyes}\binits{F.}} \AND
\bauthor{\bsnm{Sz-Nagy},~\bfnm{B{\'e}la}\binits{B.}}
(\byear{1990}).
\btitle{Functional Analysis}.
\bpublisher{Dover}, \baddress{New York}.
%Reprint of
% the 1955 original}.
\bid{mr={1068530}}
\bptok{imsref}%
\end{bbook}
%
\endbibitem

%b28 #&#
\bibitem{row00}
%
\begin{barticle}[pbm]
\bauthor{\bsnm{Roweis},~\bfnm{S.~T.}\binits{S.~T.}} \AND
\bauthor{\bsnm{Saul},~\bfnm{L.~K.}\binits{L.~K.}}
(\byear{2000}).
\btitle{Nonlinear dimensionality reduction by locally linear embedding}.
\bjournal{Science}
\bvolume{290}
\bpages{2323--2326}.
\bid{doi={10.1126/science.290.5500.2323}, issn={0036-8075},
pii={290/5500/2323}, pmid={11125150}}
\bptok{imsref}%
\end{barticle}
%
\endbibitem

%b29 #&#
\bibitem{spe98}
%
\begin{barticle}[pbm]
\bauthor{\bsnm{Spellman},~\bfnm{P.~T.}\binits{P.~T.}},
\bauthor{\bsnm{Sherlock},~\bfnm{G.}\binits{G.}} \AND
\bauthor{\bsnm{Zhang},~\bfnm{M.~Q.}\binits{M.~Q.}}
(\byear{1998}).
\btitle{Comprehensive identification of cell cycle-regulated genes of
the yeast
Saccharomyces cerevisiae by microarray hybridization}.
\bjournal{Mol. Biol. Cell}
\bvolume{9}
\bpages{3273--3297}.
\bid{issn={1059-1524}, pmcid={25624}, pmid={9843569}}
\bptok{imsref}%
\end{barticle}
%
\endbibitem

%b30 #&#
\bibitem{tan09}
%
\begin{barticle}[pbm]
\bauthor{\bsnm{Tang},~\bfnm{Rong}\binits{R.}} \AND
\bauthor{\bsnm{M{\"{u}}ller},~\bfnm{Hans-Georg}\binits{H.-G.}}
(\byear{2009}).
\btitle{Time-synchronized clustering of gene expression trajectories}.
\bjournal{Biostatistics}
\bvolume{10}
\bpages{32--45}.
\bid{doi={10.1093/biostatistics/kxn011}, issn={1468-4357}, pii={kxn011},
pmid={18502728}}
\bptok{imsref}%
\end{barticle}
%
\endbibitem

%b31 #&#
\bibitem{ten00}
%
\begin{barticle}[pbm]
\bauthor{\bsnm{Tenenbaum},~\bfnm{J.~B.}\binits{J.~B.}}, \bauthor
{\bparticle{de}
\bsnm{Silva},~\bfnm{V.}\binits{V.}} \AND
\bauthor{\bsnm{Langford},~\bfnm{J.~C.}\binits{J.~C.}}
(\byear{2000}).
\btitle{A global geometric framework for nonlinear dimensionality reduction}.
\bjournal{Science}
\bvolume{290}
\bpages{2319--2323}.
\bid{doi={10.1126/science.290.5500.2319}, issn={0036-8075},
pii={290/5500/2319}, pmid={11125149}}
\bptok{imsref}%
\end{barticle}
%
\endbibitem

%b32 #&#
\bibitem{tud54}
%
\begin{barticle}[auto:STB|2011/12/15|13:36:40]
\bauthor{\bsnm{Tuddenham},~\bfnm{R.}\binits{R.}} \AND
\bauthor{\bsnm{Snyder},~\bfnm{M.}\binits{M.}}
(\byear{1954}).
\btitle{Physical growth of California boys and girls from
birth to
age 18}.
\bjournal{California Publications on Child Development}
\bvolume{1}
\bpages{183--364}.
\bptok{imsref}%
\end{barticle}
%
\endbibitem

%b33 #&#
\bibitem{wan99}
%
\begin{barticle}[mr]
\bauthor{\bsnm{Wang},~\bfnm{Kongming}\binits{K.}} \AND
\bauthor{\bsnm{Gasser},~\bfnm{Theo}\binits{T.}}
(\byear{1999}).
\btitle{Synchronizing sample curves nonparametrically}.
\bjournal{Ann. Statist.}
\bvolume{27}
\bpages{439--460}.
\bid{doi={10.1214/aos/1018031210}, issn={0090-5364}, mr={1714722}}
\bptok{imsref}%
\end{barticle}
%
\endbibitem

%b34 #&#
\bibitem{yao05}
%
\begin{barticle}[mr]
\bauthor{\bsnm{Yao},~\bfnm{Fang}\binits{F.}},
\bauthor{\bsnm{M{\"u}ller},~\bfnm{Hans-Georg}\binits{H.-G.}} \AND
\bauthor{\bsnm{Wang},~\bfnm{Jane-Ling}\binits{J.-L.}}
(\byear{2005}).
\btitle{Functional data analysis for sparse longitudinal data}.
\bjournal{J. Amer. Statist. Assoc.}
\bvolume{100}
\bpages{577--590}.
\bid{doi={10.1198/016214504000001745}, issn={0162-1459}, mr={2160561}}
\bptok{imsref}%
\end{barticle}
%
\endbibitem

\end{thebibliography}
\end{document}